\documentclass[a4paper,12pt]{amsart}
\usepackage{xcolor}
\usepackage{fullpage}
\usepackage{amsmath,amssymb,amsthm,amsfonts}
\usepackage{graphicx}
\usepackage{enumerate}
\usepackage{bbm}
\usepackage{hyperref}
\usepackage{float}
\restylefloat{table}
\usepackage{multirow}
\usepackage{subfig}
\usepackage{placeins}
\usepackage{enumitem}

\numberwithin{equation}{section}

\newtheorem{algorithm}{Algorithm}[section]

\renewcommand{\div}{\operatorname{div}}

\def\b1{{\mathbf 1}}
\def\bv{{\mathbf v}}
\def\bu{{\mathbf u}}

\def\bg{{\mathbf g}}

\def\br{{\mathbf r}}

\def\bx{{\mathbf x}}
\def\by{{\mathbf y}}

\def\bbf{{\mathbf f}}

\def\bF{{\mathbf F}}

\newcommand{\N}{{\mathcal N}}

\newcommand{\T}{{\mathcal T}}

\renewcommand{\div}{\text{div}}

\newcommand{\norm}[1]{\left\|#1\right\|}
\title[DNN enhanced nonlinear discretizations and solvers]
{DNN{\tiny} Approximation of Nonlinear Finite Element Equations}
\author{Tuyen Tran$^2$, Aidan Hamilton, Maricela Best McKay, Benjamin Quiring, and Panayot S. Vassilevski$^{1,2}$}

\address{${}^1$Center for Applied Scientific Computing,
             Lawrence Livermore National Laboratory,
             P.O. Box 808, L-561,
            Livermore, CA 94551, U.S.A.}

\email{tuyen2@pdx.edu, hamilton49@llnl.gov, bestmckay1@llnl.gov, quiring1@llnl.gov, vassilevski1@llnl.gov}

\address{${}^2$Fariborz Maseeh Department of Mathematics and Statistics, Portland State University, Portland, Oregon, USA}

\email{tuyen2@pdx.edu, panayot@pdx.edu}

\thanks{This work was performed under the auspices of the U.S. Department of Energy by Lawrence Livermore National Laboratory under Contract DE-AC52-07NA27344.}

\begin{document}
 
\begin{abstract}
We investigate the potential of applying (D)NN ((deep) neural networks) for approximating nonlinear mappings 
arising in the finite element discretization of nonlinear PDEs (partial differential equations).
As an application, we apply the trained DNN to replace the coarse nonlinear operator
thus avoiding the need to visit the 
fine level discretization in order to evaluate the actions of the true coarse nonlinear operator. The feasibility of 
the studied approach is demonstrated in a two-level FAS (full approximation scheme) used to solve a nonlinear diffusion-reaction PDE.
\end{abstract}

\maketitle

\section{Introduction}

In recent times deep neural networks, (\cite{Goodfellow}), have become the method of choice in solving state of the art machine learning problems, such as classification, clustering, pattern recognition, and prediction with enormous impact in many applied areas.  
There is also an increasing trend in scientific computing to take advantage of the potential of DNNs as nonlinear approximation tool, (\cite{ApproxNNs}). This goes both for using DNNs in devising new approximation algorithms as well as for trying to develop  mathematical theories that explain and quantify the ability of DNNs as universal approximation methodology, with results originating in (\cite{Cybenko1989}) to many recent works, especially in the area of convolutional NN, (\cite{siegel2019approximation}, \cite{zhou2018universality},  \cite{daubechies2019nonlinear}). 
At any rate, this is still a very active area of research with no ultimate theoretical result available yet. 

Recently, the deep neural networks have been also utilized in the field of numerical solution of PDEs (\cite{bar2019}, \cite{MM}, \cite{HanJE17}, \cite{hsieh2019learning}, \cite{ke2016multigrid}), and for convolutional ones, see (\cite{He2019}). 

In this work, we investigate the ability of fully connected DNNs to provide an accurate enough approximation of the nonlinear mappings that arise in the finite element discretization of nonlinear PDEs. The finite element method applied to a nonlinear PDE posed variationally, basically requires the evaluation of integrals over each finite element which involves nonlinear functions that can be evaluated pointwise (generally, using quadratures). The unknown function $u_h$ which approximates the solution of the PDE  is a linear combination of piecewise polynomials and the individual integrals for any given value of $u_h$ can be 
evaluated accurately enough (in general, approximately, using quadrature formulas). 
Typically, in a finite element discretization procedure, we use refinement. That is, we can have a fine enough final mesh, a set of elements $\T_h$, obtained from a previous level coarse mesh $\T_H$. 
One way to solve the fine-level nonlinear discretization problem, is to utilize the  existing hierarchy of discretizations. One approach to maintain high accuracy of the coarse operators while evaluating their actions is to utilize the accurate fine level nonlinear operator. That is, for any given coarse finite element function, we expand it in terms of the fine level basis, apply the action of the fine level nonlinear operator and then restrict the result back to the coarse level, i.e., we apply a Galerkin procedure. This way of defining the coarse nonlinear operator 
provides better accuracy, however its evaluation requires fine-level computations. 
In the linear case, one can actually precompute the coarse level operators (matrices) explicitly, which is not the case in the nonlinear case. This study offers a way to generate a coarse operator, that can approximate the  variationally defined finite element one on coarse levels by training a fully connected DNN. We do not do this globally, 
but rather construct the desired nonlinear mapping based on actions of locally trained DNNs associated with each coarse element from a coarse mesh $\T_H$. This is much more feasible than training the actions of the global 
coarse nonlinear mapping; these will have as input many coarse functions (i.e., their coefficient vectors) and are thus much bigger and hence much more expensive than their restrictions to the individual coarse elements.

The remainder of this paper is structured as follows.
In Section~\ref{section: problem setting in general terms}, we introduce the problem in a general setting. Then in the following Section~\ref{subsection: training results}, we present  a computational study of the training of the coarse nonlinear operators by varying the domain (a box in a high-dimensional space with size equal to the number of local coarse degrees of freedom). The purpose of the study is to  assess the complexity of the local DNNs
depending on the desired approximation accuracy.
We also show the approximation accuracy of the global coarse nonlinear mapping (of our main interest) which depends on the coarsening ratio $H/h$. 
In Section~\ref{section: FAS}, we introduce the FAS (full approximation scheme) solver (\cite{Brandt}), and in the following section~\ref{subsection: FAS results}, we apply the trained DNNs to replace the true coarse operator in a two-level FAS  for a model nonlinear diffusion-reaction PDE discretized by piecewise linear elements.
Finally, in Section~\ref{conclusion}, we draw some conclusions and outline few directions for possible future work.

\section{Approximation for nonlinear mappings using DNNs} \label{section: problem setting in general terms}

\subsection{Problem setting}\label{subsection: problem setting}

We are given the system of nonlinear equations 
\begin{equation}\label{nonlinear problem}
F(\bu) = \bbf.
\end{equation}
Here, $F$ is a mapping from ${\mathbb R}^n \mapsto {\mathbb R}^n$, and we have access only at its actions (by calling some function). 

We assume that the solution belongs to a box $K \subset {\mathbb R}^n$, e.g.,
$\bu \in [-a,a]^n$ for some value of $a >0$. Typically, nonlinear problems like \eqref{nonlinear problem} are solved by iterations, and for any given current iterate $\bu$, we look for a correction $\bg$ such that 
$\bu:=\bu+\bg$ gives a better accuracy. This motivates us to rewrite \eqref{nonlinear problem}
 as
\begin{equation*}
G(\bu,\;\bg) = {\overline \bbf},
\end{equation*}
where
\begin{equation*}
G(\bu,\;\bg):= F(\bu+\bg) - F(\bu) \text{ and }{\overline \bbf}:= \bbf - F(\bu).
\end{equation*}
Our goal is to train a DNN where $\bu$ and $\bg$ are the inputs and  $G(\bu,\; \bg)$ is the output.
The input $\bu$ is drawn from the box $K$, whereas the correction $\bg$ is drawn from a small
ball $B=\{\|\bg\| \le \delta\}$. In our study to follow, we vary the parameters $a$ and $\delta$ for a particular mapping $F$ (and respective $G$) to assess the complexity of the resulting DNN and examine the approximation accuracy. The general strategy is as follows. We draw $m_a \ge 1$ vectors from the box $K$ using Sobol sequence (\cite{Sobol,sobol2}) and $m_\delta \ge 1$ vectors from the ball $B$, also using Sobol sequence. The alternative would be to simply use random points in $K$ and $B$, however Sobol sequence is better in terms of cost versus approximation ability (at least for smooth mappings).
Once we have built the DNN with a desired accuracy on the training data, we test its approximation quality on a number of randomly selected points from $K$ and $B$.

Our results are documented in the next subsection for a particular example of a finite element mapping; first on each individual subdomain (coarse element) $T \in \T_H$ and then for its global action composed from all locally trained DNNs.

\subsection{Training DNNs for model nonlinear finite element mapping}\label{subsection: results for training coarse nonlinear mappings}

We consider the nonlinear PDE
\begin{equation}\label{Poisson}
\begin{array}{rl}
-\div (k(u) \nabla u)+u&=f\,\mbox{on}\,\Omega,\\
\nabla u \cdot \vec{n} &=0\, \mbox{on}\,\partial\Omega.
\end{array}
\end{equation}
Here, $\Omega$ is a polygon in ${\mathbb R}^2$ and $k(u)$ is a given positive nonlinear function of $u$. 

The variational formulation for \eqref{Poisson} is: find $u \in H^1(\Omega)$ such that 
\begin{equation*}
\int_{\Omega}(k(u)\nabla u \cdot\nabla v +uv)\;dx = \int_{\Omega} fv\;dx\, \mbox{for all} \,v \in H^1(\Omega).
\end{equation*}
The above problem is discretized by piecewise linear finite elements on triangular mesh $\T_h$ that yields a system of nonlinear equations.

In this section, we consider the coarse nonlinear mapping that corresponds to a coarse triangulation $\T_H$ which after refinement gives the fine one $\T_h$. The coarse finite element space is $V_H$ and the fine one is $V_h$.
By construction, we have $V_H \subset V_h$. Let $\{\phi^H_i\}^{N_H}_{i=1}$ be the basis of $V_H$ and 
$\{\phi^h_i\}^{N_h}_{i=1}$ be the basis of $V_h$. These are piecewise linear functions associated with their respective 
triangulations $\T_H$ and $\T_h$. More specifically, we use  Lagrangian bases, i.e., $\phi^H_i$ and 
$\phi^h_i$ are associated with the sets of vertices, $\N_H$ and $\N_h$, of the elements of their respective triangulations. 

The coarse nonlinear operator $F:=F_H$ is then defined as follows.
Let $u_H\in V_H$  be a coarse finite element function. Since $V_H \subset V_h$, we can expand $u_H$ in terms of the basis of $V_h$, i.e.,
\begin{equation*}
u_H= \sum\limits_{\bx_i \in \N_h} u_H(\bx_i) \phi^h_i.
\end{equation*}
 We can also expand $u_H$ in terms of the basis of $V_H$, i.e., we have
\begin{equation*}
u_H  = \sum\limits_{\bx_i \in \N_H} u_H(\bx_i) \phi^H_i.
\end{equation*}
In the actual computations we use their coefficient vectors
\begin{equation*}
\bu_c = (u_H(\bx_i))_{\bx_i \in \N_H} \in {\mathbb R}^{N_H} \text{ and }
\bu = (u_H(\bx_i))_{\bx_i \in \N_h} \in {\mathbb R}^{N_h}.
\end{equation*}
These coefficient vectors are related by an interpolation mapping $P$ (which is piecewise linear),
i.e., we have 
\begin{equation*}
\bu = P \bu_c.
\end{equation*}

First we define the local nonlinear mappings $F:= F^H_T$, associated with each $T \in \T_H$. 

In terms of finite element functions, we have as input
$u_H$ restricted to $T$, and we evaluate the integrals
\begin{equation*}
F^H_T:\;\left . u_H \right |_T \mapsto \int\limits_T k(u_H) \nabla u_H \cdot\nabla \phi^H_i\;d\bx,
\text{ for all vertices } \bx_i \text { of the coarse element } T.
\end{equation*}
Each integral over $T$ is computed as a sum of integrals over the fine-level elements 
$\tau \subset T$, $\tau \in \T_h$, using fine level computations, i.e., $u_H$ and $\phi^H_i$  are linear on 
each $\tau$ and these fine-level integrals are assumed computable (by the finite element software used to generate the fine level discretization, which  possibly employs high order quadratures).

In terms of linear algebra computations, we have $\bu_{c,\;T} = (u^H(\bx_i))_{\bx_i \in \N_H \cap T}$
 as an input vector, and have as an output
a vector $F^H_T(\bu_{c,\;T})$ of the same size, i.e., equal to  the number of vertices of $T$. 
Note that we will be training the DNN for the mapping of two variables, $\bu_{c,\;T}$ and $\bg_{c,\:T}$, i.e.,
\begin{equation*}
G_T(\bu_{c,\;T},\;\bg_{c,\;T}):= F^H_T(\bu_{c,\;T}+\bg_{c,\;T}) - F^H_T(\bu_{c,\;T}).
\end{equation*}
 That is, the input vectors will have size two times bigger than the output 
vectors.
Once we have trained the local actions of the nonlinear mapping, the global action is obtained by standard {\em assembly},
using the fact that 
\begin{equation}
\label{global coarse mapping assembled from local ones}
G(\bu_c,\;\bg_c) = \sum\limits_{T\in \T_H} I_T G_T(\bu_{c,\;T},\;\bg_{c,\;T}).
\end{equation}
Here, $I_T$ stands for the  mapping that extends a local vector defined on $T$ to a global vector by zero values outside $T$. 
For a global vector $\bv_c$, $\bv_{c,\;T} = (I_T)^T \bv_c = \left . \bv_c \right |_T$ denotes its restriction to $T$.

In the following section, we provide actual test results for training DNNs, first for the local mappings $G_T$, and then for the respective global one $G$.

\subsection{Training local DNNs for the model finite element mapping}
\label{subsection: training results}

We use \textbf{Keras}, (\cite{pal2017deep}), a Python Deep Learning library to approximate nonlinear mappings. \textbf{Keras }is a high-level neural networks API (Application Programming Interface), written in Python and capable of running on top of TensorFlow (\cite{shukla2018machine}). In this work, we use a fully connected network. The \textit{Sequential model}  in \textbf{Keras} provides a way to implement such a network.

The input vectors are of size $2n_c$  ($\bu_{c,T}$ stacked on top of $\bg_{c,T}$) and the desired outputs are the actions ($G_T(\bu_{c,T},\;\bg_{c,T})$) represented as vectors of size $n_c$ for any given input. The network consists of few fully connected layers with \textbf{tanh} activation at each layer. We use the standard mean squared error as the loss function.

In the tests to follow, we use data $\bu_{c,\;T}$ from the boxes $K$, and $\bg_{c,T}$ from  balls 
$B$ of various sizes.
Specifically, we performed the numerical tests with $10$ $\bu_{c,\;T}$ vectors each taken from 
$$K=[-1,1]^{n_c}, [-0.1,0.1]^{n_c}, [-0.05,0.05]^{n_c}, [-0.01,0.01]^{n_c},$$
and $50$ $\bg_{c,\;T}$ vectors drawn from balls $B$ with radii $\delta_B = 0.1, 0.05, 0.01,0.005 $, respectively (i.e., the first $K$ is paired with the first ball $B$, the second box $K$ is paired with the second ball $B$, and so on).  
In our test we have chosen $n_c=4$, that is, the local sets $T$ have four coarse dofs. 
Also, we vary the ratio $H/h = 2,4,8$, which implies that we have $9,16,$ and $81$ fine dofs, respectively (while keeping fixed the number of coarse dofs $n_c =4$ in $T$).  

The network was trained with $3$ layers, each with $16$ neurons. The training algorithm was provided by TensorFlow using the ADAM optimizer (\cite{Kingma2014AdamAM}) which is a variant of the SGD (stochastic gradient descent) algorithm, see, e.g., (\cite{Gilbert}). 
We used $500$ epochs, $batch\,size=10$, learning rate $\alpha=0.001$ along with $\beta_1 =0.9$ and $\beta_2=0.999$. For more details on the meaning of these parameters, we refer to (\cite{pal2017deep}) and (\cite{Gilbert}).

\noindent\begin{minipage}{\linewidth}
\centering
\captionsetup{type=figure}
\includegraphics[scale=0.5]{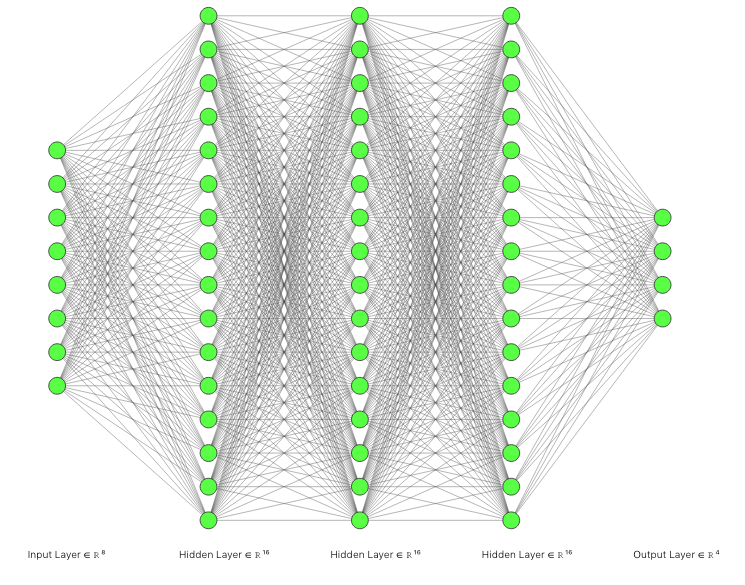} 
\captionsetup{type=figure}
\captionof{figure}{Schematic representation of the used network architecture}
\label{NNs}
\end{minipage}

\vspace{1cm}

Let $G_{T_{DNN}}$ be the action after the training. Tables \ref{table H/h=2}, \ref{table H/h=4}, and \ref{table H/h=8} show  the average of the relative errors using $\norm{}_2$ and $\norm{}_\infty$ of $\dfrac{\norm{G_T(\bu_{c,\;T},\;\bg_{c,\;T})-G_{T_{DNN}}(\bu_{c,\;T},\;\bg_{c,\;T})}_2}{\norm{G_T(\bu_{c,\;T},\;\bg_{c,\;T})}_2}$ and 
$\dfrac{\norm{G_T(\bu_{c,\;T},\;\bg_{c,\;T})-G_{T_{DNN}}(\bu_{c,\;T},\;\bg_{c,\;T})}_{\infty}}{\norm{G_T(\bu_{c,\;T},\;\bg_{c,\;T})}_{\infty}}$ over $100$ examples consisting of $10$ $\bu_{c,\;T}$ within the box $K$ and $10$ $\bg_{c,\;T}$ within the ball $B$.
\vspace{1cm}

\noindent\begin{minipage}{\linewidth}
\centering
\begin{tabular}{ |c|c|c|c|}
 \hline
$K$ & $\delta_B$ & relative $\ell_2$-error& relative $\ell_{\infty}$-error\\
  \hline
$[-1,1]^{n_c}$ & 0.1 & 0.005209&0.007234 \\
  \hline 
  $[-1,1]^{n_c}$ & 0.05 & 0.000829&0.001168  \\
  \hline 
  $[-1,1]^{n_c}$ & 0.01 & 0.000245&0.000374  \\
  \hline 
  $[-1,1]^{n_c}$ & 0.005 & 0.000125
&0.000193  \\
  \hline 
$[-0.1,0.1]^{n_c}$ & 0.05 & 0.000131&0.000186 \\ 
 \hline
 $[-0.1,0.1]^{n_c}$ & 0.01 & 3.858919E-06& 5.973702E-06\\ 
 \hline
 $[-0.1,0.1]^{n_c}$ & 0.005 &3.966236E-06 &6.276684E-06 \\ 
 \hline
 $[-0.05,0.05]^{n_c}$ & 0.01 &1.043632E-06&1.487640E-06  \\ 
 \hline
 $[-0.05,0.05]^{n_c}$ & 0.005 &6.841004E-07 &1.048721E-06  \\ 
 \hline
  $[-0.01,0.01]^{n_c}$ & 0.005 & 3.636231E-08&5.151947E-08\\ 
 \hline
\end{tabular}
\captionsetup{type=table}
\captionof{table}{The relative average $L_2$ and $L_{\infty}$ errors for $H/h=2$}
\label{table H/h=2}
\end{minipage}

\noindent\begin{minipage}{\linewidth}
\centering
\begin{tabular}{ |c|c|c|c|}
 \hline
$K$ & $\delta_B$ & relative $\ell_2$-error& relative $\ell_{\infty}$-error\\
  \hline
$[-1,1]^{n_c}$ & 0.1 & 0.001198&0.001670 \\
  \hline 
  $[-1,1]^{n_c}$ & 0.05 & 0.0007022&0.001021 \\
  \hline 
  $[-1,1]^{n_c}$ & 0.01 & 0.000254&0.000377  \\
  \hline 
  $[-1,1]^{n_c}$ & 0.005 & 1.692172E-05
&2.455399E-05  \\
  \hline 
$[-0.1,0.1]^{n_c}$ & 0.05 & 1.692172E-05&2.455399E-05\\ 
 \hline
 $[-0.1,0.1]^{n_c}$ & 0.01 & 3.105803E-06& 4.837250E-06\\ 
 $[-0.1,0.1]^{n_c}$ & 0.005 & 2.251858E-06&3.573355E-06 \\ 
 \hline
 $[-0.05,0.05]^{n_c}$ & 0.01 & 8.565291E-07& 1.322711E-06 \\ 
 \hline
 $[-0.05,0.05]^{n_c}$ & 0.005 & 6.394531E-07& 9.717701E-07 \\ 
 \hline
  $[-0.01,0.01]^{n_c}$ & 0.005 & 2.872041E-08&4.041187E-08\\ 
 \hline
\end{tabular}
\captionsetup{type=table}
\captionof{table}{The relative average $L_2$ and $L_{\infty}$ errors for $H/h=4$}
\label{table H/h=4}
\end{minipage}

\noindent\begin{minipage}{\linewidth}
\centering
\begin{tabular}{ |c|c|c|c|}
 \hline
$K$ & $\delta_B$ & relative $\ell_2$-error& relative $\ell_{\infty}$-error\\
  \hline
$[-1,1]^{n_c}$ & 0.1 &0.001060& 0.001010 \\
  \hline 
  $[-1,1]^{n_c}$ & 0.05 &0.000641 &0.001003\\
  \hline 
  $[-1,1]^{n_c}$ & 0.01 &  0.000196& 0.000292\\
  \hline 
  $[-1,1]^{n_c}$ & 0.005 & 1.141443E-05
& 1.723346E-05 \\
  \hline 
$[-0.1,0.1]^{n_c}$ & 0.05 &1.151866E-05 &1.744020E-05\\ 
 \hline
 $[-0.1,0.1]^{n_c}$ & 0.01 &2.756513E-06 &4.175218E-06\\ 
 \hline
 $[-0.1,0.1]^{n_c}$ & 0.005 & 2.115907E-06 & 3.441013E-06\\ 
 \hline
 $[-0.05,0.05]^{n_c}$ & 0.01 &8.378942E-07& 1.221053E-06 \\ 
 \hline
 $[-0.05,0.05]^{n_c}$ & 0.005 &5.626164E-07 &8.913452E-07 \\ 
 \hline
  $[-0.01,0.01]^{n_c}$ & 0.005 &2.645869E-08 &3.999300E-08\\ 
 \hline
\end{tabular}
\captionsetup{type=table}
\captionof{table}{The relative average $L_2$ and $L_{\infty}$ errors for $H/h=8$}
\label{table H/h=8}
\end{minipage}

The approximation of the global coarse nonlinear operator is presented in Table~\ref{table: global coarse nonlinear operator approximation}.
As it is seen from formula \eqref{global coarse mapping assembled from local ones}, we can approximate each $G_T$ independently of each other, and after combining the individual approximations (using the same assembly formula with each $G_T$ replaced by its DNN approximation), we define the approximation to the global $G$.
We use the same setting for training the individual neural networks for each $G_T$. We have decomposed $\Omega$ into several subdomains $T \in \T_H$ so that which each $T$ has $n_c=4$. In this test,  we chose $H/h=2$. 

In Table \ref{table: global coarse nonlinear operator approximation}, we show  how the accuracy varies for different local boxes $K$ and respective balls $B$.
As before, we present the average of the relative $L_2$ and $L_{\infty}$ errors.
One can notice that for finer $h$ (and $H = 2h$), i.e., more local subdomains $T$, 
we get somewhat better approximations of the global $G$. It should be expected since \
with smaller $h$ the finite element problem approximates better the continuous one. 
Some visual illustration of this fact is presented in Figures \ref{figure: 4 subdomains, min L_2 error},~\ref{figure: 4 subdomains, max L_2 error},~\ref{figure: 16 subdomains, min L_2 error},~\ref{figure: 16 subdomains, max L_2 error},~\ref{figure: 64 subdomains, min L_2 error}, and \ref{figure: 64 subdomains, max L_2 error}.
More specifically, we provide plots of $G_T$ and $G_{T_{DNN}}$ for the data that achieves the min and max $L_2$ errors when the number of subdomains are  4, 16 and 64 corersponding to box $K=[-0.05,0.05]$ and and ball $B$ with $\delta_B=0.005$.

\vspace{1 cm}

\noindent\begin{minipage}{\linewidth}
\centering
\begin{tabular}{ |c|c|c|c|} 
 \hline
\multirow{ 2}{*}{\# of subdomains} & $K=[-1,1]$&$K=[-0.1,0.1]$&$K=[-0.05,0.05]$\\
&$\delta_B=0.1$&$\delta_B=0.01$&$\delta_B=0.005$\\
  \hline
  \multirow{ 2}{*}{4}& $L_2\,\mbox{error:}\,$0.004838&6.379939E-05&9.578247E-06\\
  &$L_{\infty}\,\mbox{error:}\,$ 0.008154&0.000130&3.154789E-05\\
 \hline
   \multirow{ 2}{*}{16}& $L_2\,\mbox{error:}\,$0.003847&4.363349E-05&1.112007E-06\\
  &$L_{\infty}\,\mbox{error:}\,$ 0.005229&0.000021&1.828478E-05\\
  \hline
    \multirow{ 2}{*}{64}& $L_2\,\mbox{error:}\,$0.002977&7.388419E-06&8.997734E-07\\
  &$L_{\infty}\,\mbox{error:}\,$ 0.003401&2.789513e-5&3.907758E-06\\
  \hline
    \multirow{ 2}{*}{100}& $L_2\,\mbox{error:}\,$: 
    0.000909&5.655502e-6&9.043519E-08\\
  &$L_{\infty}\,\mbox{error:}\,$ 0.001688&1.542110e-5&4.098461E-07\\
  \hline
\end{tabular}
\captionsetup{type=table}
\captionof{table}{Comparison for different numbers of subdomains\label{table: global coarse nonlinear operator approximation}}
\end{minipage}

For the next test, we have $n_c=4$ and $H/h=4$. The same settings for neural networks are used along with $K=[-0.05,0.05]$ and $\delta_B=0.005$. The results in Table \ref{table: H/h=4} show the average of the relative $L_2$ and $L_{\infty}$ errors for different number of subdomains.
\vspace{1 cm} 

\noindent\begin{minipage}{\linewidth}
\centering
\begin{tabular}{|c|c|c|c|c|}
\hline
\# of subdomains & $L_2$ min        & $L_2$ max   & $L_2$           &$L_{\infty}$   \\ \hline
4         & 2.204530E-06  & 1.012479E-04 & 5.276626E-06 & 2.367189E-05     \\ \hline
16        & 1.307212E-07  & 2.768105E-05 & 7.107769E-07 & 2.068912E-06     \\ \hline
64        & 4.8551896E-08 & 1.395091E-05 & 5.341155E-08 & 7.789968E-07 \\ \hline
\end{tabular}
\captionsetup{type=table}
\captionof{table}{The relative average $L_2$ and $L_{\infty}$ errors for $H/h=4$}
\label{table: H/h=4}
\end{minipage}

As expected, since smaller $h$ and fixed $H$ gives better approximation of the nonlinear operator, for the ratio $H/h=4$ we do see better approximation than for $H/h=2$.

\begin{figure}[!ht]
\centering
\subfloat[$G_T(\bu_{c,\;T},\;\bg_{c,\;T})$]{\includegraphics[scale=0.25]{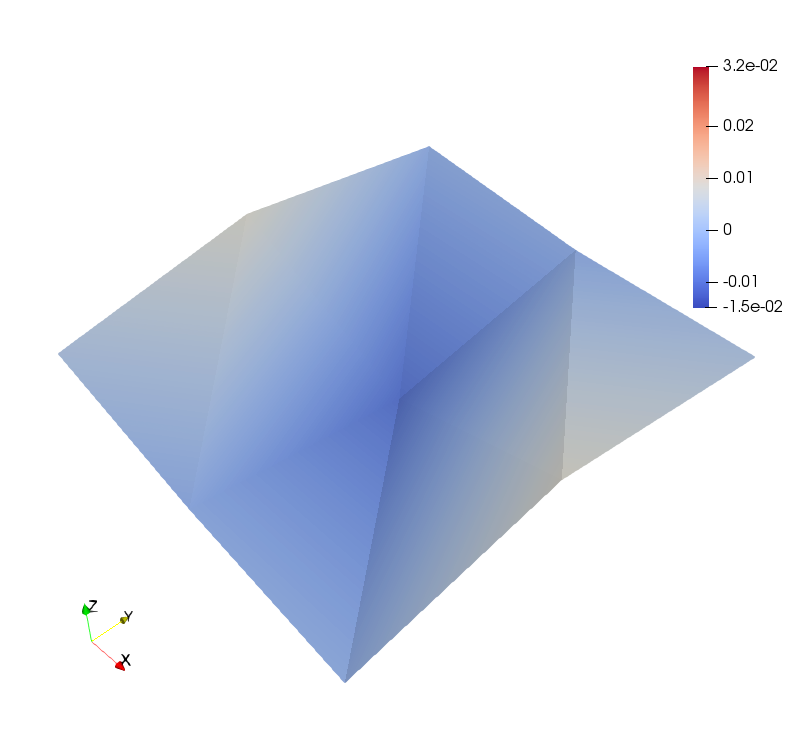}} 
\subfloat[$G_{T_{DNN}}(\bu_{c,\;T},\;\bg_{c,\;T})$]{\includegraphics[scale=0.25]{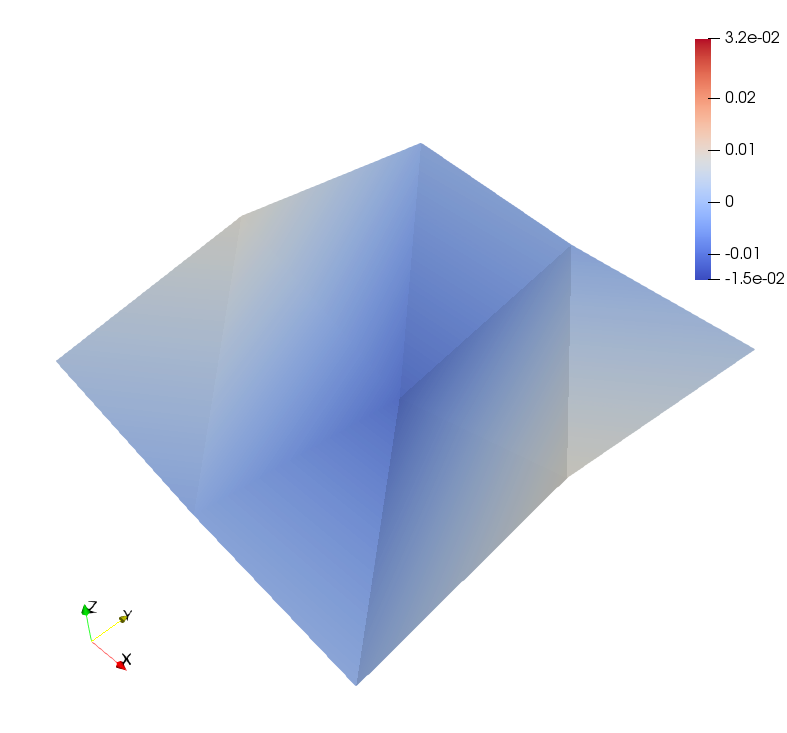}}
\caption{Plots of $G_T$ and $G_{T_{DNN},\;\bg_{c,\;T}}$ for four subdomains at the $(\bu_{c,\;T},\;\bg_{c,\;T})$ which achieves the min $L_2$ error. The min relative error is $2.662757E{-06}$.}
\label{figure: 4 subdomains, min L_2 error}
\end{figure}

\newpage
\begin{figure}[!ht]
\centering
\subfloat[$G_T(\bu_{c,\;T},\;\bg_{c,\;T})$]{\includegraphics[scale=0.25]{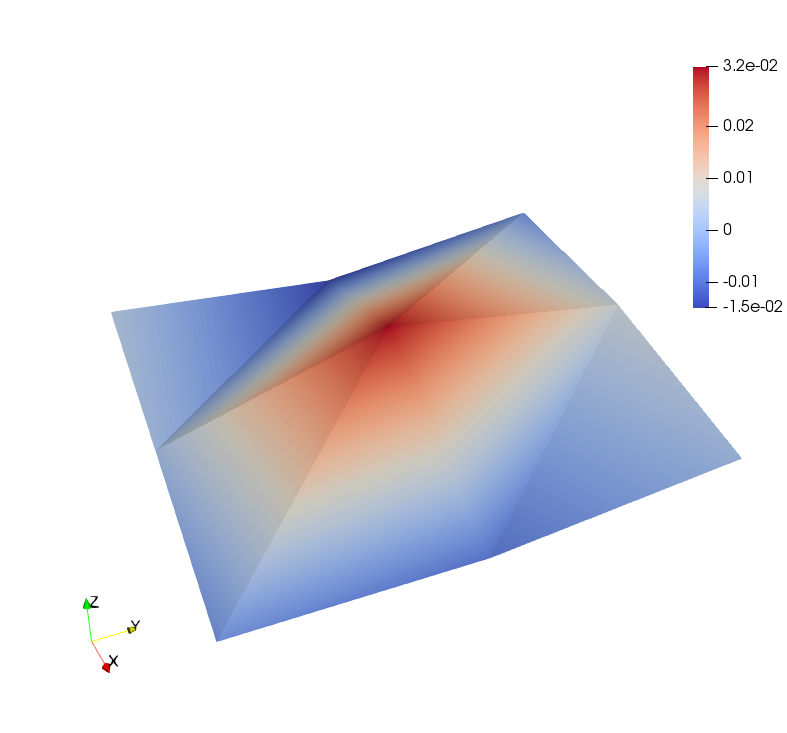}} 
\subfloat[$G_{T_{DNN}}(\bu_{c,\;T},\;\bg_{c,\;T})$]{\includegraphics[scale=0.25]{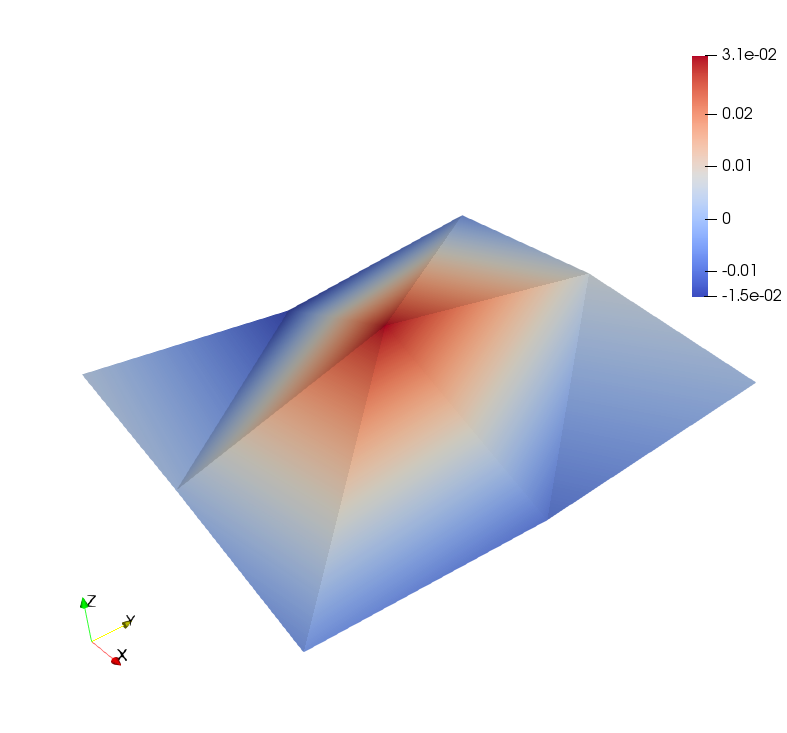}}
\caption{Plots of $G_T$ and $G_{T_{DNN}}$ for four subdomains at the $(\bu_{c,\;T},\;\bg_{c,\;T})$ which achieves the max $L_2$ error.  The max relative error is $1.122974E{-04}$.}
\label{figure: 4 subdomains, max L_2 error}
\end{figure}

\begin{figure}[!ht]
\centering
\subfloat[$G_T(\bu_{c,\;T},\;\bg_{c,\;T})$]{\includegraphics[scale=0.25]{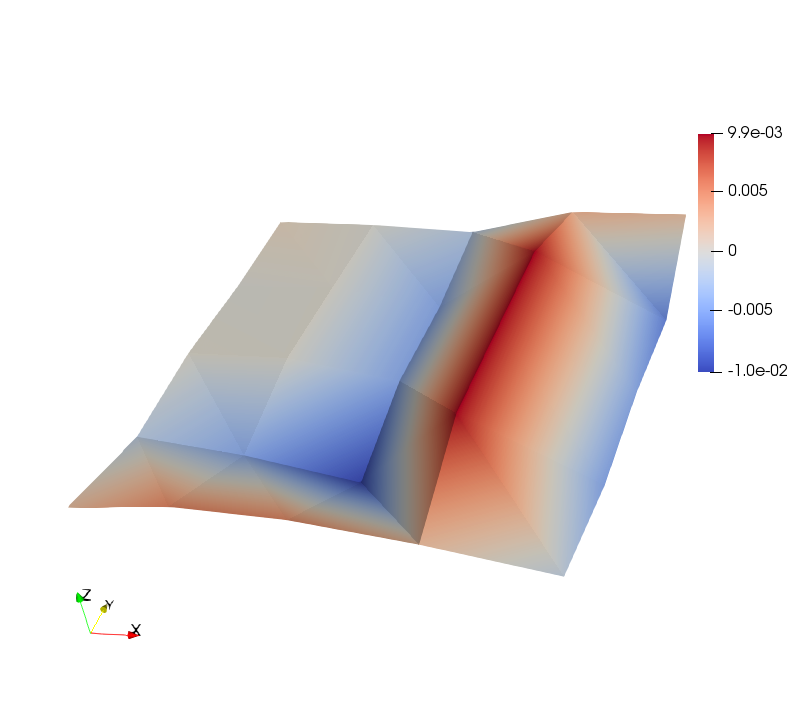}} 
\subfloat[$G_{T_{DNN}}(\bu_{c,\;T},\;\bg_{c,\;T})$]{\includegraphics[scale=0.25]{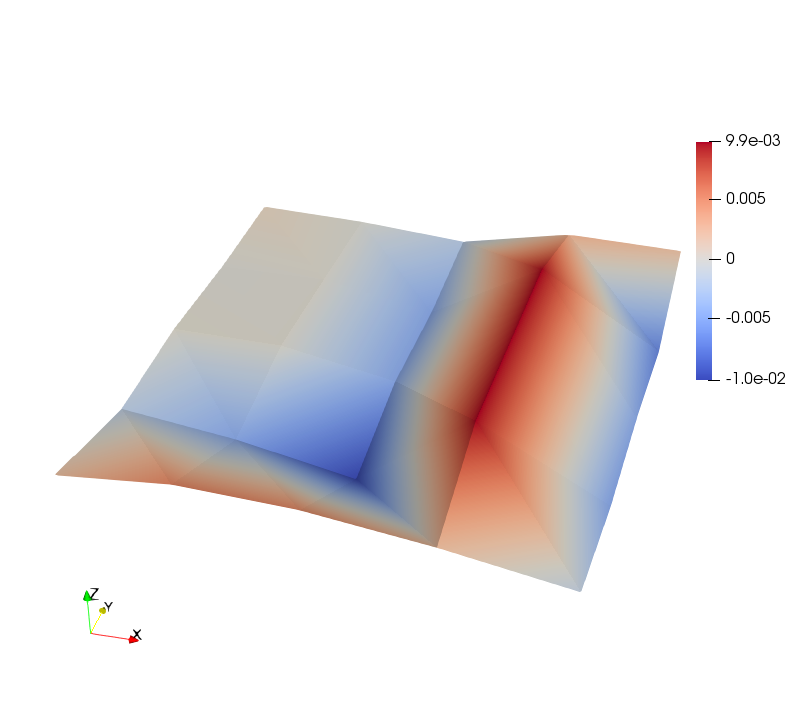}}
\caption{Plots of $G_T$ and $G_{T_{DNN}}$ for sixteen subdomains at the $(\bu_{c,\;T},\;\bg_{c,\;T})$ which achieves the min $L_2$ error. The min relative error is $8.886265E{-07}$. }
\label{figure: 16 subdomains, min L_2 error}
\end{figure}

\newpage
\begin{figure}[!ht]
\centering
\subfloat[$G_T(\bu_{c,\;T},\;\bg_{c,\;T})$]{\includegraphics[scale=0.3]{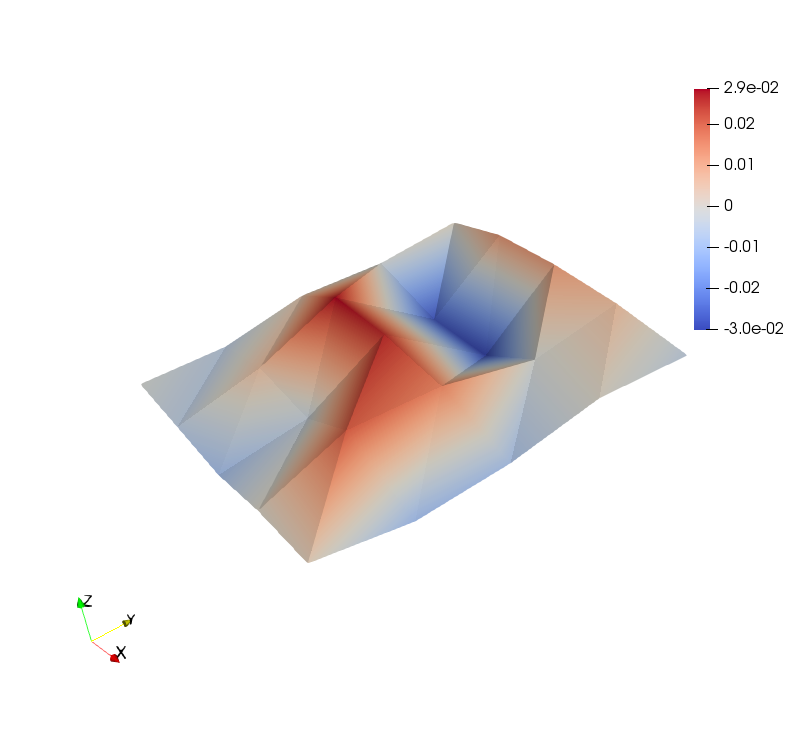}} 
\subfloat[$G_{T_{DNN}}(\bu_{c,\;T},\;\bg_{c,\;T})$]{\includegraphics[scale=0.3]{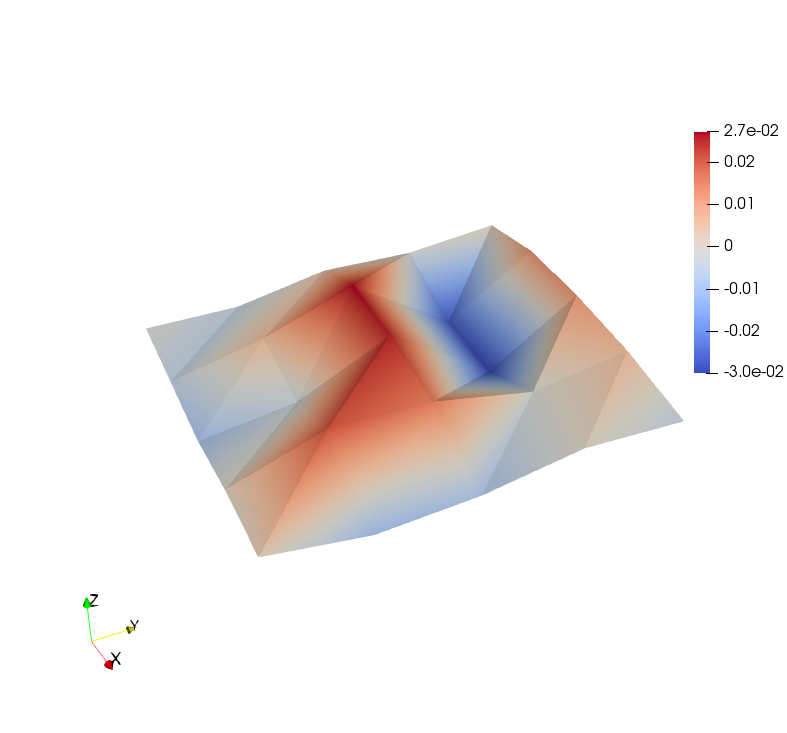}}
\caption{Plots of $G_T$ and $G_{T_{DNN}}$ for sixteen subdomains at the $(\bu_{c,\;T},\;\bg_{c,\;T})$ which achieves the max $L_2$ error. The max relative error is $1.285776E{-05}$.}
\label{figure: 16 subdomains, max L_2 error}
\end{figure}

\newpage
\begin{figure}[!ht]
\centering
\subfloat[$G_T(\bu_{c,\;T},\;\bg_{c,\;T})$]{\includegraphics[scale=0.12]{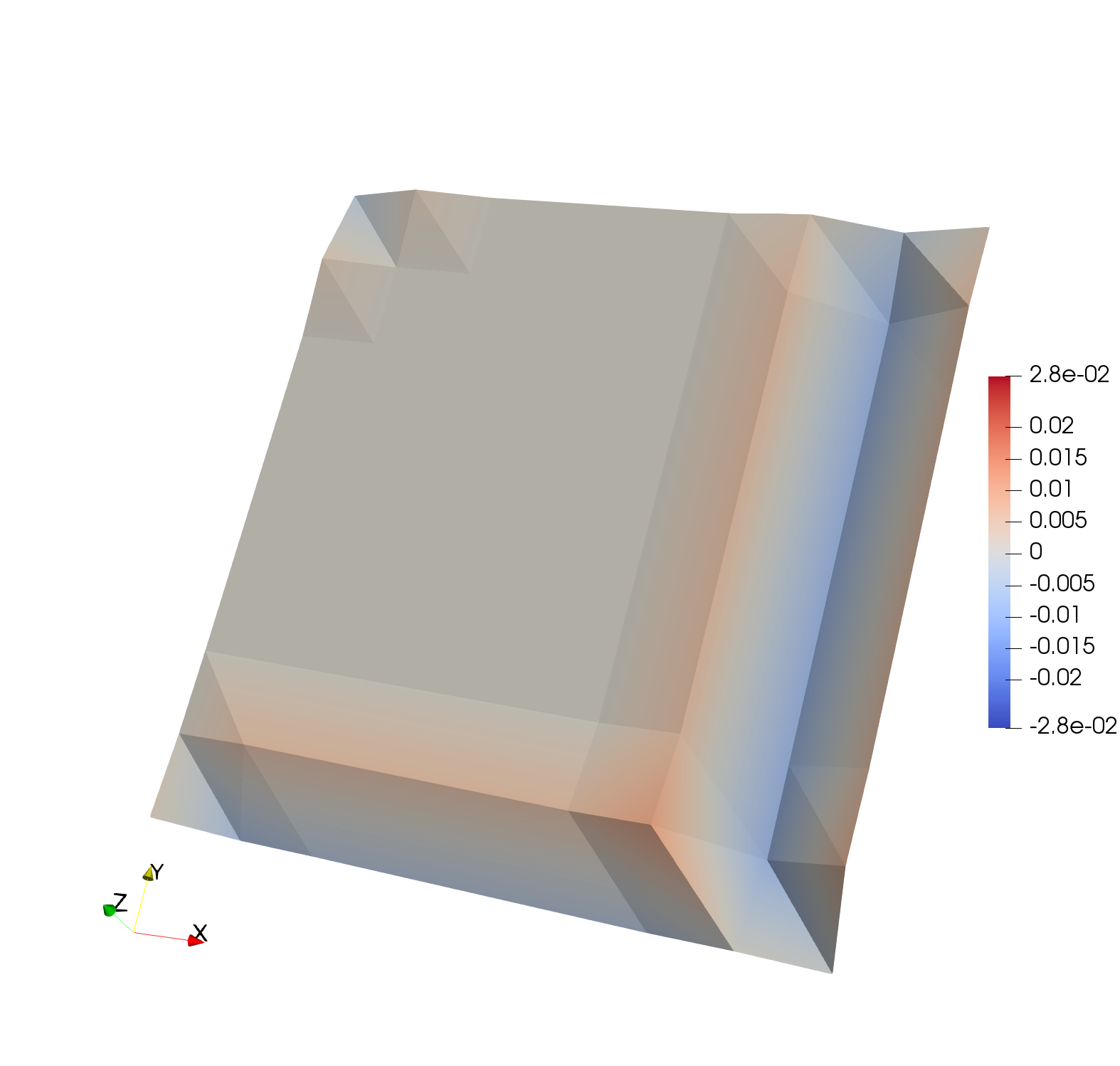}} 
\subfloat[$G_{T_{DNN}}(\bu_{c,\;T},\;\bg_{c,\;T})$]{\includegraphics[scale=0.12]{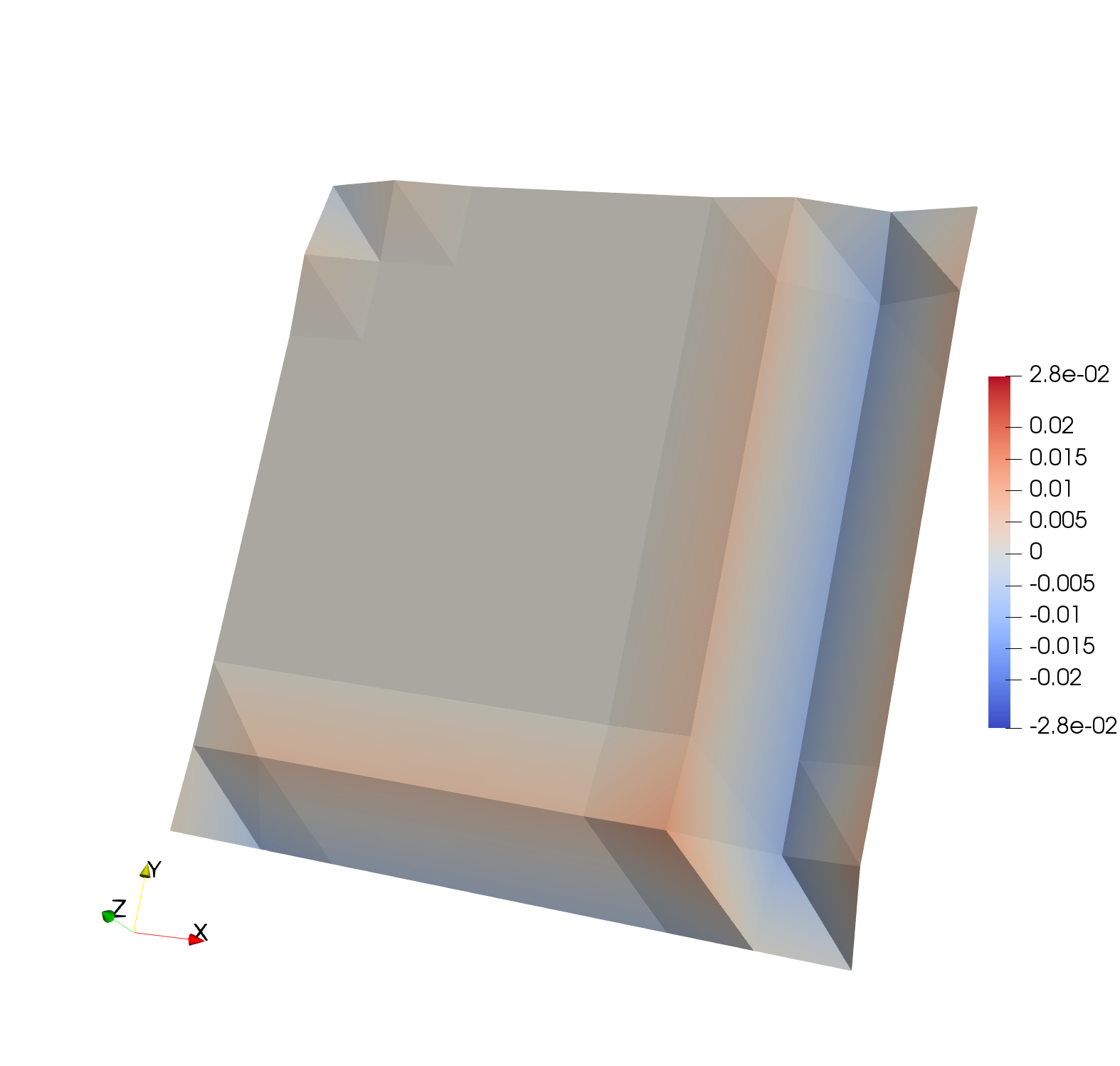}}
\caption{Plots of $G_T$ and $G_{T_{DNN}}$ for sixty four subdomains at the $(\bu_{c,\;T},\;\bg_{c,\;T})$ which achieves the min $L_2$ error. The min relative error is $1.493536E{-07}$.}
\label{figure: 64 subdomains, min L_2 error}
\end{figure}

\begin{figure}[!ht]
\centering
\subfloat[$G_T(\bu_{c,\;T},\;\bg_{c,\;T})$]{\includegraphics[scale=0.12]{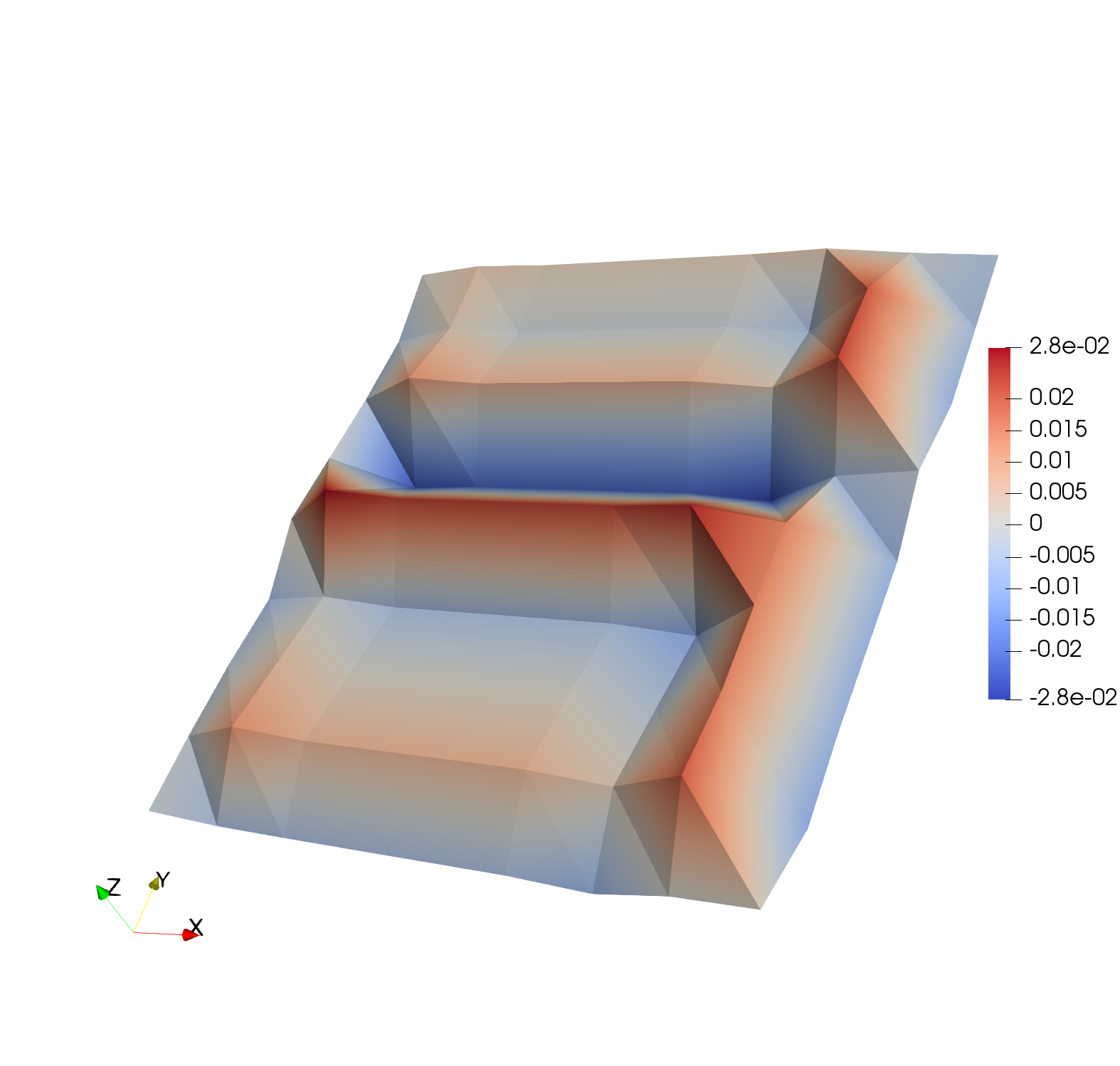}} 
\subfloat[$G_{T_{DNN}}(\bu_{c,\;T},\;\bg_{c,\;T})$]{\includegraphics[scale=0.12]{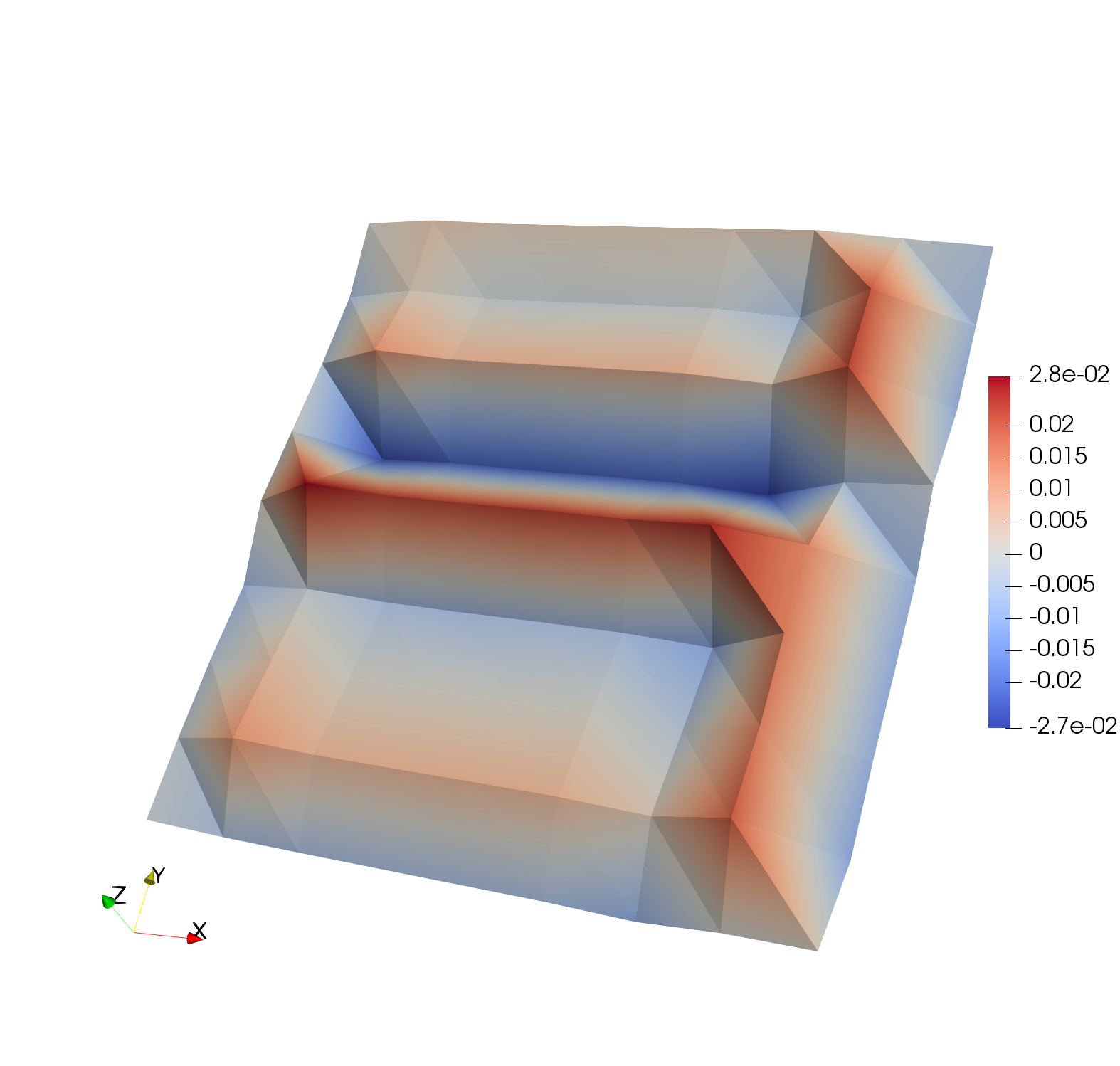}}
\caption{Plots of $G_T$ and $G_{T_{DNN}}$ for sixty four subdomains at the $(\bu_{c,\;T},\;\bg_{c,\;T})$ which achieves the max $L_2$ error. The max relative error is $2.717447E{-05}$.}
\label{figure: 64 subdomains, max L_2 error}
\end{figure}

\newpage
\section{Application of DNN approximate coarse mappings in two-level FAS}
\label{section: FAS}

\subsection{The FAS algorithm}\label{FAS-algorithm}
A standard approach for solving \eqref{nonlinear problem} is to use Newton's method. The latter is an iterative process in which given a current iterate $\bu$, we compute the next one, $\bu_{next}$, by solving the Jacobian equation
\begin{equation}\label{Jacobian residual equation}
J_F(\bu) \by = \br:=\bbf - F(\bu),
\end{equation}
and then $\bu_{next} = \bu + \by$.

Typically, the Jacobian problem
\eqref{Jacobian residual equation} is solved by an iterative method such as GMRES (generalized minimal residual).

To speed-up the convergence, for nonlinear equations coming from finite element discretizations of elliptic PDEs, such as \eqref{Poisson}, we can exploit hierarchy of discretizations.
A popular method is the two-level FAS (full approximation scheme) proposed by Achi Brandt, (\cite{Brandt}) (see also (\cite{vassilevski})). For a recent convergence analysis of FAS, we refer to (\cite{chen2018convergence}). 

To define the two-level FAS, we need a coarse version of $F$, which is another nonlinear mapping $F_c:\; {\mathbb R}^{n_c} \mapsto {\mathbb R}^{n_c}$, for some $n_c < n$. We also need its Jacobian $J_c = J_{F_c}$.
Again, we assume that both are available via their actions (by calling some appropriate functions). To communicate data between the original, {\em fine level}, ${\mathbb R}^n$ and the {\em coarse level}, ${\mathbb R}^{n_c}$, we are given two linear mappings (matrices):
\begin{itemize}
\item [(i)] Coarse-to-fine (interpolation or prolongation) mapping $P:\; {\mathbb R}^{n_c} \mapsto {\mathbb R}^n$.
\item [(ii)] Fine-to-coarse projection $\pi:\; {\mathbb R}^n \mapsto {\mathbb R}^{n_c}$, more 
precisely, we assume that $\pi P = I$ (then $(P\pi)^2 = P\pi$ is a projection).
In our finite element setting, $\pi$ is simply the restriction of a fine-grid vector $\bu$ to the coarse 
dofs, i.e., $\pi = [0,\;I]$, where the columns in $I$ correspond to the coarse dofs viewed as subset of the the fine dofs. 
\end{itemize}

Then, the two-level FAS (TL-FAS) can be formulated as follows.
\begin{algorithm}[Two-level FAS]\label{algorithm: TL FAS} \hfill

For problem \eqref{nonlinear problem}, with a current approximation $\bu$, the 
two-level FAS method  performs the following steps to compute the next approximation 
$\bu_{next}$.
\begin{itemize}

\item For a given $m \ge 1$ apply $m$ steps of (inexact) Newton algorithm, \eqref{Jacobian residual equation}, to compute 
$\by_m$ and let $\bu_{\frac{1}{3}} = \bu + \by_m$.
\item Form the coarse-level nonlinear problem for $\bu_c$
\begin{equation}\label{coarse nonlinear problem}
F_c(\bu_c) = \bbf_c \equiv F_c(\pi \bu_{\frac{1}{3}}) + P^T(\bbf - F(\bu_{\frac{1}{3}})).
\end{equation}
\item Solve \eqref{coarse nonlinear problem} accurately enough using Newton's method 
based on the coarse Jacobian $J_c$ and initial iterate $\bu_c:= \pi \bu_{\frac{1}{3}}$.
 Here we use enough iterations of GMRES for solving the 
resulting coarse Jacobian residual equations 
\begin{equation*}
J_c(\bu_c) \by_c = \br_c = \bbf_c - F_c(\bu_c), \text{ and let } \bu_c:= \bu_c + \by_c,
\end{equation*}
until a desired accuracy is reached. 
\item Update fine-level approximation
\begin{equation*}
\bu_{\frac{2}{3}} = \bu_{\frac{1}{3}}+P (\bu_c-\pi \bu_{\frac{1}{3}}).
\end{equation*}
\item Repeat the FAS cycle starting with $\bu:=\bu_{next} = \bu_{\frac{2}{3}}$ 
until a desired accuracy is reached.
\end{itemize}
\end{algorithm}

In what follows, we use the following equivalent form of the TL-FAS. At the coarse level, we will represent
$\bu_c := \bu^0_c + \bg_c$, where $\bu^0_c = \pi \bu_{\frac{1}{3}}$ is the initial coarse iterate coming from the fine level, and we will be solving for the correction $\bg_c$. That is, the coarse problem in terms of the correction reads
\begin{equation*}
{\overline F}_c(\bg_c)  = {\overline \bbf}_c,
\end{equation*}
where
\begin{equation*}
\begin{array}{rl}
{\overline F}_c(\bg_c) &\equiv G_c(\bu^0_c,\;\bg_c) =
F_c(\bu^0_c+ \bg_c) - F_c(\bu^0_c) \\
{\overline \bbf}_c &\equiv P^T(\bbf - F(\bu_{\frac{1}{3}})) = \bbf_c - F_c(\bu^0_c).
\end{array}
\end{equation*}
The rest of the algorithm does not change, in particular, we have
\begin{equation*}
\begin{array}{rl}
\br_c &= \bbf_c - F_c(\bu_c) \\
& = \bbf_c - F_c(\bu^0_c+\bg_c)\\
& = {\overline \bbf}_c + F_c(\bu^0_c) - F_c(\bu^0_c+\bg_c)\\
& = {\overline \bbf}_c - G_c(\bu^0_c,\;\bg_c)\\
&= {\overline \bbf}_c - {\overline F}_c(\bg_c),
\end{array}
\end{equation*}
and
\begin{equation*}
\bu_{\frac{2}{3}} = \bu_{\frac{1}{3}} + P (\bu_c- \pi \bu_{\frac{1}{3}}) = 
\bu_{\frac{1}{3}} + P(\bu_c-\bu^0_c) =
\bu_{\frac{1}{3}}  + P \bg_c.
\end{equation*}

\subsubsection{The choice for $F_c$ using DNN}

In our finite element setting, a true (Galerkin) coarse operator is $P^TF(P(.))$, where $P$ is the piecewise linear interpolation mapping and $F$ is the fine level nonlinear finite element operator.
We train the global coarse actions based on the fact that the actions of 
$F$ and $P$ can be computed subdomain-by-subdomain employing standard finite element assembly procedure, as described in the previous section. That is, $F$ can be assembed by local 
$F_T$s and the coarse $P^TF(P(.))$ can be assembled from local coarse actions $P^T_T F_T (P_T(.))$ based on local versions $\{P_T\}$ of $P$.

More specifically, we train for each subdomain $T \in \T_H$ a DNN which takes as input any pair of  coarse vectors $\bv_{c,T},\;\bg_{c,\;T}  \in {\mathbb R}^{n_c}$ and produces $P^T_TF_T(P_T(\bv_{c,T}+\bg_{c,T}))-P^T_TF_T(P_T\bv_{c,T}) \in {\mathbb R}^{n_c}$ as our desired output. 
The global action $P^TF(P(\bv_c+\bg_c))- P^TF(P(\bv_c))$ is computed by assembling all local actions, and we use the same assembling procdeure for the approximations obtained using the trained local DNNs.

The trained this way DNN gives the actions of our coarse nonlinear mapping ${\overline F}_c(.)$.


\subsection{Some details on implementing TL-FAS algorithm}

We are solving the nonlinear equation $F(\bu) = \bbf$ where the actions of $F:\; {\mathbb R}^n \mapsto {\mathbb R}^n$  
are available. Also, we assume that we can compute its Jacobian matrix for any given $\bu$, $J(\bu)$. 

We are also given an interpolation matrix $P:\;{\mathbb R}^{n_c} \mapsto {\mathbb R}^n$. 
Finally, we have a mapping $\pi:\; {\mathbb R}^n \mapsto {\mathbb R}^{n_c}$
such that $\pi P = I$ on ${\mathbb R}^{n_c}$. This implies that 
$P \pi$ is a projection, i.e, $(P\pi)^2 = P\pi$.

Consider the coarse nonlinear mapping $F(\bu_c) \equiv P^T(F(P\bu_c)):\;
{\mathbb R}^{n_c} \mapsto {\mathbb R}^{n_c}$. 

We assume that an approximation $G(\bv_c,\;\bg_c)$ to the mapping 
\begin{equation*}
F_c(\bv_c+\bg_c) - F_c(\bv_c),
\end{equation*}
is available through its actions for a set of input vectors $\bv_c$ varying in a given box $K$  and for another input vector $\bg_c$ varying in  a small ball $B$ about the origin. For a fixed $\bv_c$, we denote ${\overline F}_c(\bg_c) =
G(\bv_c,\;\bg_c)$. 

We are interested in the following two-level FAS algorithm for solving 
$F(\bu) =\bbf$.

\begin{algorithm}[TL-FAS] \label{algorithm: TL-FAS-2}
\hfill

Input parameters:

\begin{itemize}
\item Initial approximation $\bu_0$ sufficiently close to the exact solution $\bu_*$. For a problem with a known solution, we choose $\bu_0 = \bu_* + \tau\times( \text{random vector})$, where $\tau$ is an input (e.g., $\tau =1,\; 0.1,\; 10.0, ...$). The random vector has as components random numbers in $(-1,1)$.

\item $\delta$ (e.g., $\delta = 0.1$, or $\delta = 0.5$) - tolerance used in GMRES to solve approximately the fine-level Jacobian equations.

\item Maximal number $N_{\max} = 1,2 \text{ or } 4,$ of fine-level inexact Newton iterations

\item Additionally, maximal number of GMRES iterations, $I_{\max} = 2, \text{ or }5,$  allowed in solving the  fine-level Jacobian equations.

\item $\delta_c$ (e.g., $\delta_c = 10^{-3}$), tolerance used in GMRES for solving coarse-level Jacobian equations. 

\item $\tau_c$ (equal to $1$ or $0.1$, or $0.01$) - step length in coarse-level inexact Newton iterations.

\item Maximal number of coarse-level GMRES iterations, $I^c_{\max} = 1000$. 

\item Maximal number of inexact coarse-level Newton iterations ,
$N^c_{\max} = 10 \text{ or } 100$.

\item Maximal number of FAS iterations, $N_{FAS} = 10$.
\item Tolerance for FAS iterations, $\epsilon = 10^{-6}$. 
\end{itemize}
With the above input specified, a TL-FAS algorithm takes the form:

\begin{itemize}
\item FAS-loop: If visited $< N_{FAS}$ times perform the steps below. Otherwise 
exit.

\begin{itemize}
\item 
Perform $N_{\max}$ fine-level inexact Newton iterations, which involve the following steps: 

\begin{itemize}
\item For the current iterate $\bu$ (the initial one, $\bu = \bu_0$,  is given as input)
compute residual 
\begin{equation*}
\br = \bbf - F(\bu).
\end{equation*}
\item Compute Jacobian $J(\bu)$.
\item Solve approximately the Jacobian equation
\begin{equation*}
J(\bu) \by = \br,
\end{equation*}
using at most $I_{\max}$ GMRES iterations or until reaching tolerance $\delta$, i.e., 
\begin{equation*}
\|J(\bu) \by - \br\| \le \delta \; \|\br\|.
\end{equation*}
\item Update $\bu:= \bu + \by$.
\end{itemize}
\end{itemize}

\begin{itemize}
\item Compute fine-level Jacobian $J(\bu)$ and the 
coarse-level one  $J_c = P^TJ(\bu)P$. 
\item Compute $\bu_c= \pi \bu$.
\item Coarse-loop for solving
\begin{equation*}
{\overline \bF}_c(\bg_c):= G(\bu_c,\bg_c) = {\overline \bbf}_c \equiv P^T (\bbf - F(\bu)),
\end{equation*}
with initial guess $\bg_c=0$ where we keep $\bu_c$  and the coarse Jacobian 
$J_c$ fixed. The coarse-level loop reads:
\setlength{\abovedisplayskip}{3pt}
\setlength{\belowdisplayskip}{3pt} 
\begin{itemize}
\item Compute coarse residual
\begin{equation*}
\br_c = {\overline \bbf}_c - {\overline F}_c(\bg_c).
\end{equation*}
\item Solve by GMRES, $J_c \by_c = \br_c$ using at most $I^c_{\max}$ iterations 
or until we reach $\|J_c \by_c -\br_c\| \le \delta_c \;\|\br_c\|.$
\item Update
\begin{align*}
\bg_c:= \bg_c + \tau_c \by_c.
\end{align*}
\item Repeat at most $N^c_{\max}$ times the above three steps of the coarse-level loop.
\end{itemize}
\item Update fine-level iterate
\begin{equation*}
\bu:=\bu + P \bg_c.
\end{equation*}
\item If $\|F(\bu) - \bbf\| > \epsilon \|F(\bu_0) - \bbf)\|$, go to the beginning of  FAS-loop.
Otherwise exit.

\end{itemize}
\end{itemize}

\end{algorithm}

\subsection{Local tools for FAS}

\begin{figure}[!ht]
\centering
\subfloat[Global coarse mesh]{\includegraphics[scale=0.15]{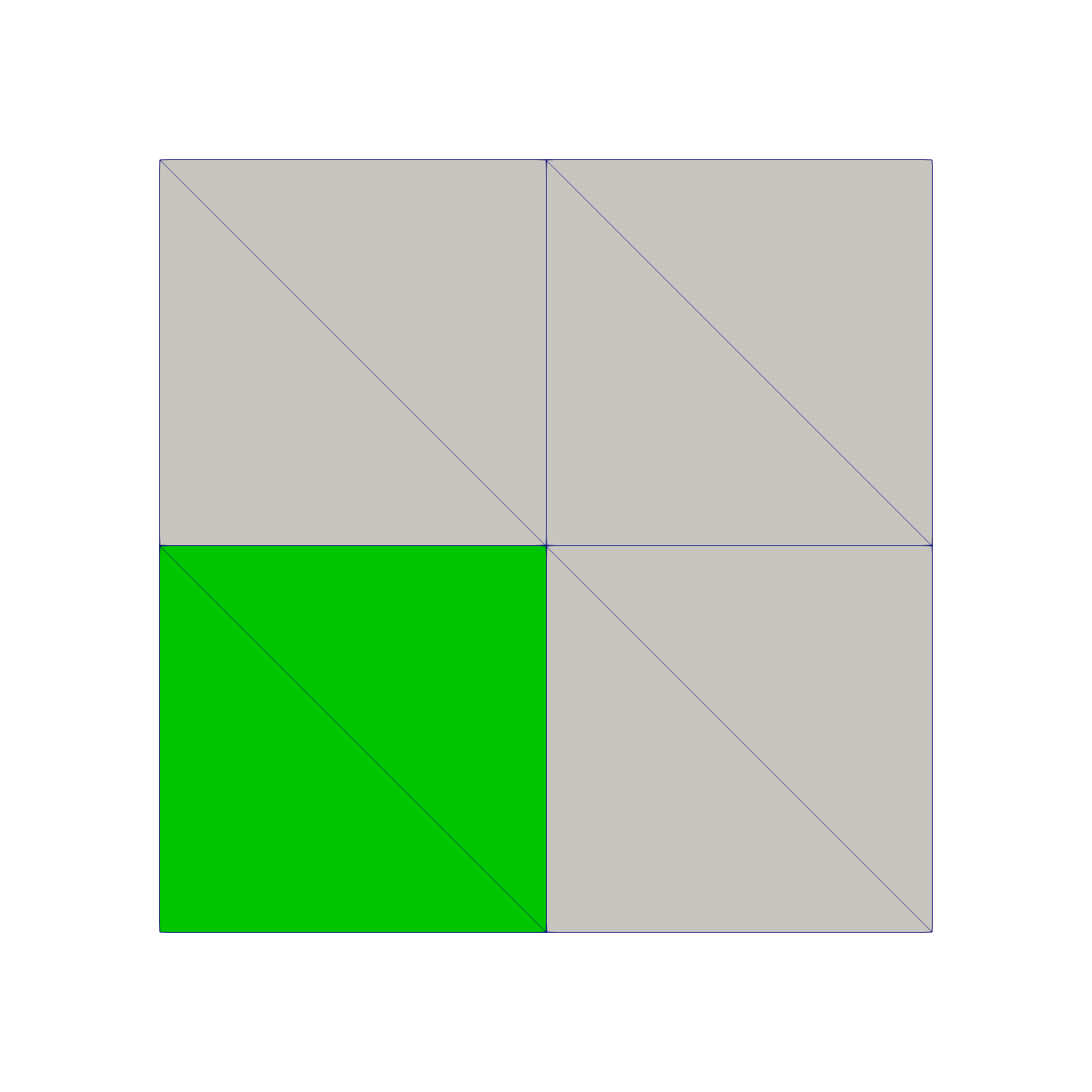}} 
\subfloat[Global fine mesh]{\includegraphics[scale=0.15]{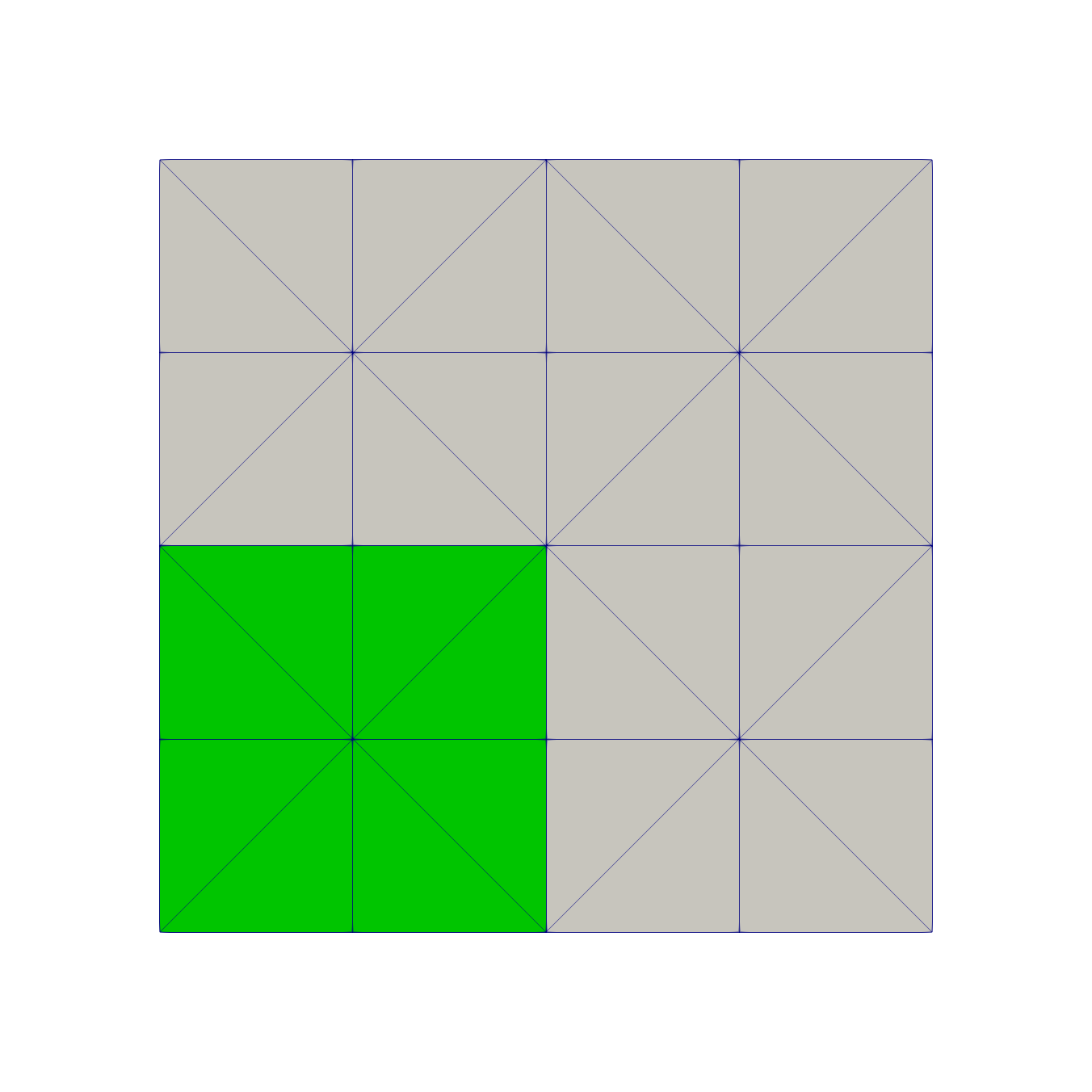}}
\caption{An example of domain decomposition for 4 subdomains when we have $H/h=2$. }
\label{mesh_2}
\end{figure}

\begin{figure}[!ht]
\centering
\subfloat[Global coarse mesh]{\includegraphics[scale=0.15]{4subdomainscoarse.png}} 
\subfloat[Global fine mesh]{\includegraphics[scale=0.15]{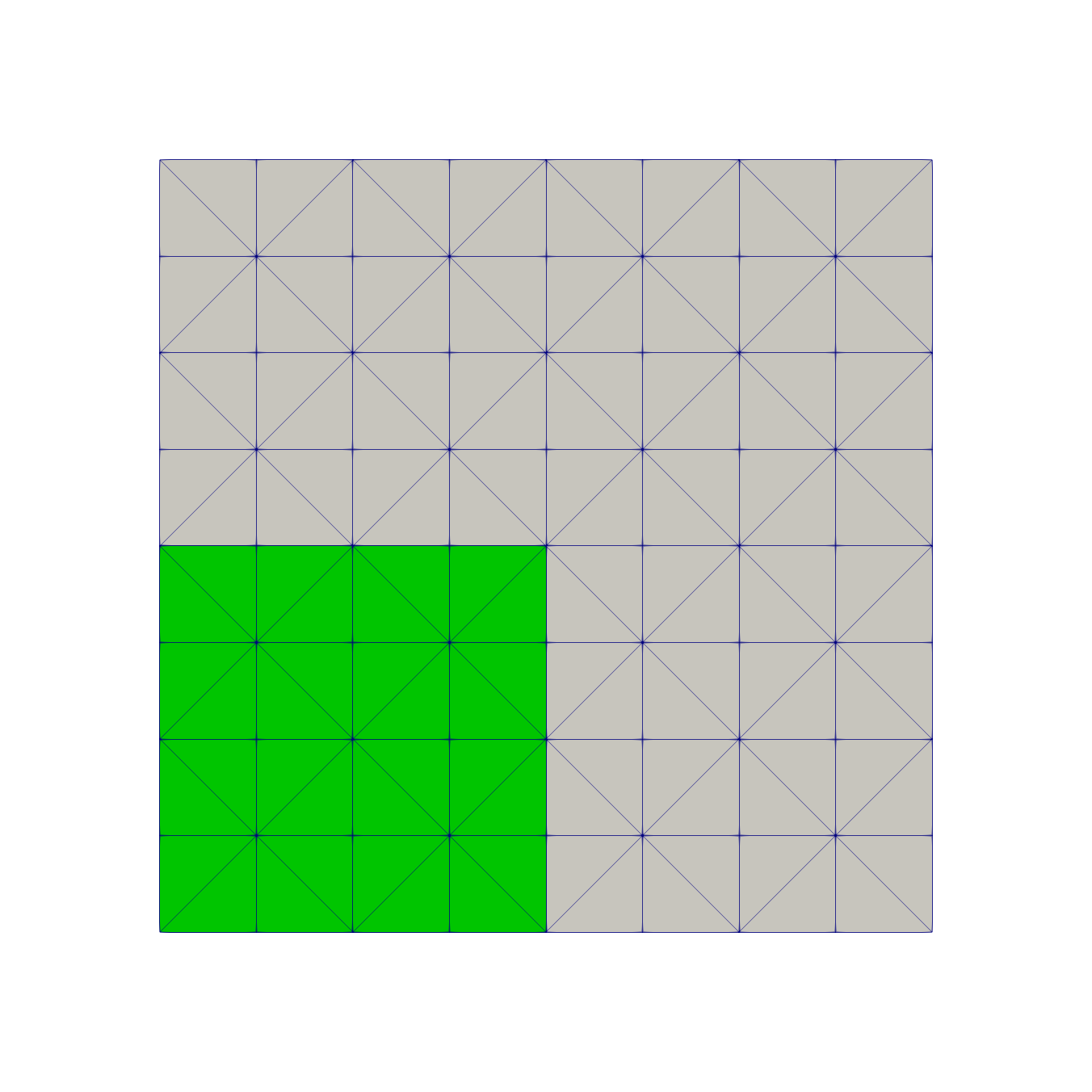}}
\caption{A visualization when we have 4 subdomains and $H/h=4$. }
\label{mesh_4}
\end{figure}

\setlength{\abovedisplayskip}{3pt}
\setlength{\belowdisplayskip}{3pt} 
We stress again that all global actions of the coarse operator (exact and approimate via DNNs) are realized by assembly of local actions. All this is possible due to the inherent nature of the finite element method. We illustrate the local subdomains in Fig.~\ref{mesh_2} and Fig.~\ref{mesh_4}.

More specifically, in Figures \ref{mesh_2} and \ref{mesh_4} ,  on the left we have an example of a global coarse mesh $\T_H$, with the shaded area being the submesh contained in $T \in \Omega_H$ -one of the four subdomains in $\T_H$.  On the right, we show the submesh  of the global fine mesh $\T_h$ restricted to the subdomain $T$. 

\subsection{Some details on implementing TL-FAS with DNNs}

Using the tools and algorithms from the previous sections, we first outline an implementation that involves training inside FAS and then present the algorithm for training outside the FAS. Training on the outside is useful because the training is then only part of a one time set up cost and independent of the r.h.s. The two-level training inside  FAS is  formulated as follows.

\begin{algorithm}[Training inside FAS]\label{algorithm:FAS inside} \hfill
With the inputs specified in algorithm \ref{algorithm: TL-FAS-2}, the training inside takes the following steps 

\begin{itemize}
\item  Perform $N_{\max}$ fine-level inexact Newton iterations, which involve  steps in algorithm \ref{algorithm: TL-FAS-2} and update $\bu$.
\end{itemize}

\begin{itemize}
\item Compute fine-level Jacobian $J(\bu)$ and the 
coarse-level one  $J_c = P^TJ(\bu)P$. 
\item Compute $\bu_c= \pi \bu$.
\item We generate samples which are in the neighborhood of $\bu_c$. More specifically, we use a shifted box
$K = \{\bu_c\} + [-0.05,0.05]$ and draw $10$ to $20$ vectors $\bu_c$ from $K$ and $30$ to $50$ vectors $\bg_c$ from the  ball $B$ with radius $\delta_B=0.005$. 
We use these sets as inputs for training the DNNs to get the approximations of $G_{c,\;T}$ for each subdomain $T \in \T_H$ and assemble them into a global coarse action, $G_c$.
Next, we enter the coarse-level loop as in Algorithm \ref{algorithm: TL-FAS-2} and the rest of the algorithm remains the same. 
Note that the training is done only for the first fine-level iterate $\bu$. 
\end{itemize}

\end{algorithm}

The following implementation gives the approximations of $G_{c,\;T}$ for each subdomain $T$  and can be performed  outside the Algorithm \ref{algorithm: TL-FAS-2}.

\begin{algorithm}[Training outside FAS]\label{algorithm:FAS outside} Given the inputs specified in algorithm \ref{algorithm: TL-FAS-2},in order to get the approximations of $G_{c,\;T}$  we proceed the training outside as follows

\begin{itemize}
\item We take $M$ inputs  $\bu$ from the box $K=[-0.05,0.05]$, and $m$ corrections $\bg$ are drawn from the ball $B$ of radius $\delta_B=0.005$. This gives a training set with $M \times n$ vectors.
\item We then use these vectors as inputs and train the neural networks which provide approximations of $G_{c,\;T}$ for each subdomain $T \in \T_H$. The latter after assembly give the global coarse action, $G_c$.
\item Next, we enter Algorithm \ref{algorithm: TL-FAS-2} with this approximate $G_c$ and the rest of the algorithm remains the same. 
\end{itemize}

\end{algorithm}

\subsection{Comparative results for FAS with exact and approximate coarse operators using DNNs}
\label{subsection: FAS results}
In this section, we present some results for two level FAS using the exact operators and its approximation from the training inside and outside as mentioned in the previous section.

We perform the tests with the neural networks using the same settings in section \ref{subsection: training results} and test our algorithm for problem \ref{nonlinear problem} with exact solution is $u_* = x^2(1-x)^2 + y^2(1-y)^2$ on the unit square  in $\mathbb{R}^2$ with $H/h=2$. For the FAS algorithm, we use the following parameters:
\begin{itemize}[topsep=2pt]
\item Initial approximation $\bu_0 = \bu_* + \tau\times( \text{random vector})$, where $\tau=2$ and the random vector has as components random numbers in $(-1,1)$.

\item $\delta = 10^{-6}$ - tolerance used in GMRES to approximately solve the fine-level Jacobian equations.

\item Maximal number $N_{\max} = 2$ of fine-level inexact Newton iterations and tolerance is $0.01$.

\item Maximal number of GMRES iterations, $I_{\max} = 4$  allowed in solving the  fine-level Jacobian equations.

\item $\delta_c = 10^{-8}$, tolerance used in GMRES for solving coarse-level Jacobian equations. 

\item $\tau_c= 0.1$ - step length in coarse-level inexact Newton iterations for FAS two levels with true operators, FAS training inside, and $\tau_c=0.0001$ for FAS training outside.

\item Maximal number of coarse-level GMRES iterations, $I^c_{\max} = 10$. 

\item Maximal number of inexact coarse-level Newton iterations ,
$N^c_{\max} = 5$ with tolerance $10^{-4}$.

\item Maximal number of FAS iterations, $N_{FAS} = 10$.
\item Tolerance for FAS iterations, $\epsilon = 10^{-6}$. 
\end{itemize}

We present the results in  Table~\ref{FAS_true} (TL-FAS with true coarse operator), Table~\ref{FAS_inside} (TL-FAS with training inside), and Table~\ref{FAS_outside} (TL-FAS with training outside). It can be seen that the training inside does give similar results to the true operators and even better in terms of relative residuals. For four subdomains, the training outside reached the same number of iterations in FAS as in Tables \ref{FAS_true} and  \ref{FAS_inside}. However, we did not achieve as small residuals as we could in the previous two cases. When we have more subdomains, the training outside requires more iterations on the coarse level than the other two approaches, but nevertheless it still meets the desired tolerance. 

\vspace{1 cm}

\noindent\begin{minipage}{\linewidth}
\centering
\begin{tabular}{|c|c|c|c|}
\hline
\multirow{3}{*}{ \# of subdomains} & \multicolumn{3}{c|}{True Operators} \\
\cline{2-4}
 & FAS iteration& \# coarse iterations&  Relative residuals \\
\hline
 4&1 &5 &0.174041\\ 
 \cline{2-4}
 &2&5&0.003939\\
  \cline{2-4}
 &3&5&4.377789E-06\\
  \cline{2-4}
 &4&1&1.081675E-11\\
\hline   
 16&1 &5 &  0.161382\\
 \cline{2-4}
 &2&5&0.001406\\
  \cline{2-4}
 &3&1&3.843124E-07\\
\hline
 64&1 &5 &  0.162089\\ 
 \cline{2-4}
 &2&5&0.002728\\
  \cline{2-4}
 &3&1&2.393705E-06\\
  \cline{2-4}
 &4&1&1.141114E-09\\
  \hline
\end{tabular}
\captionsetup{type=table}
\captionof{table}{Results for FAS using true operators with different number of subdomains using $k=1+u^2$.}
\label{FAS_true}
\end{minipage}

%
%
%

\noindent\begin{minipage}{\linewidth}
\centering
\begin{tabular}{|c|c|c|c|}
\hline
\multirow{3}{*}{ \# of subdomains} & \multicolumn{3}{c|}{Approximate Operators with inside training } \\
\cline{2-4}
  & FAS iteration& \# coarse iterations&  Relative residuals \\
\hline
 4&1 &5 &0.171590\\ 
 \cline{2-4}
 &2&5&0.003733\\
  \cline{2-4}
 &3&1&3.826213E-06\\
  \cline{2-4}
 &4&1&8.325530e-12\\
\hline   
 16&1 &5 &  0.160120\\
 \cline{2-4}
 &2&5&0.001286\\
  \cline{2-4}
 &3&1&2.946412E-07\\
\hline
 64&1 &5 & 0.159982\\ 
 \cline{2-4}
 &2&5&0.001521\\
  \cline{2-4}
 &3&1&4.671780E-07\\
  \hline
\end{tabular}
\captionsetup{type=table}
\captionof{table}{Results for training NNs inside with different number of subdomains using $k=1+u^2$.}
\label{FAS_inside}
\end{minipage}

\noindent\begin{minipage}{\linewidth}
\centering
\begin{tabular}{|c|c|c|c|}
\hline
\multirow{3}{*}{ \# of subdomains} & \multicolumn{3}{c|}{Approximate Operators with outside training } \\
\cline{2-4}
  & FAS iteration& \# coarse iterations&  Relative residuals \\
\hline
 4&1 &5 &0.186189\\ 
 \cline{2-4}
 &2&5&0.004746\\
  \cline{2-4}
 &3&5&5.300791E-06\\
  \cline{2-4}
 &4&5&5.623933e-10\\
\hline   
 16&1 &5 &  0.174228\\
 \cline{2-4}
 &2&5&0.002257\\
  \cline{2-4}
 &3&5&1.246687E-06\\
  \cline{2-4}
 &4&5&4.025110E-09\\
\hline
 64&1 &5 &  0.176886\\ 
 \cline{2-4}
 &2&5&0.004503\\
  \cline{2-4}
 &3&5&1.488866E-05\\
  \cline{2-4}
 &4&5&2.061581E-08\\
  \hline
\end{tabular}
\captionsetup{type=table}
\captionof{table}{Results for training NNs outside with different number of subdomains using $k=1+u^2$.}
\label{FAS_outside}
\end{minipage}

We also studied the case $H/h=4$ with the same setting as above for the neural networks, the same $k=1+u^2$ and the exact solution. The following tables (Tables \ref{FAS_true_H/h_4}, \ref{FAS_inside_H/h_4}, and \ref{FAS_outside_H/h_4}) represent the results for 4 and 16 subdomains. They show similar behaviour as in the previous case of $H/h = 2$.

\noindent\begin{minipage}{\linewidth}
\centering
\begin{tabular}{|c|c|c|c|}
\hline
\multirow{3}{*}{ \# of subdomains} & \multicolumn{3}{c|}{True Operators } \\
\cline{2-4}
 & FAS iteration& \# coarse iterations&  Relative residuals \\
\hline
 4&1 &5 &0.171241\\ 
 \cline{2-4}
 &2&5&0.002991\\
  \cline{2-4}
 &3&1&2.325478E-06\\
  \cline{2-4}
 &4&1&6.696561e-12\\
\hline   
 16&1 &5 &  0.168299\\ 
 \cline{2-4}
 &2&5&0.002804\\
  \cline{2-4}
 &3&1&4.457398E-06\\
  \cline{2-4}
 &4&1&3.394781E-08\\
\hline 
\end{tabular}
\captionsetup{type=table}
\captionof{table}{Results for FAS using true operators with $H/h=4$.}
\label{FAS_true_H/h_4}
\end{minipage}

\noindent\begin{minipage}{\linewidth}
\centering
\begin{tabular}{|c|c|c|c|}
\hline
\multirow{3}{*}{ \# of subdomains} & \multicolumn{3}{c|}{Approximate Operators with inside training } \\
\cline{2-4}
 & FAS iteration& \# coarse iterations&  Relative residuals \\
\hline
 4&1 &5 &0.172198\\ 
 \cline{2-4}
 &2&5&0.003446\\
  \cline{2-4}
 &3&1&1.292521E-06\\
  \cline{2-4}
 &4&1&6.055108e-12\\
\hline   
 16&1 &5 &  0.167623\\ 
 \cline{2-4}
 &2&5&0.002594\\
  \cline{2-4}
 &3&5&4.748181E-06\\
  \cline{2-4}
 &4&1&2.390020E-09\\
\hline 
\end{tabular}
\captionsetup{type=table}
\captionof{table}{Results for NNs training inside with $H/h=4$.}
\label{FAS_inside_H/h_4}
\end{minipage}

\noindent\begin{minipage}{\linewidth}
\centering
\begin{tabular}{|c|c|c|c|}
\hline
\multirow{3}{*}{ \# of subdomains} & \multicolumn{3}{c|}{Approximate Operators with outside training } \\
\cline{2-4}
 & FAS iteration& \# coarse iterations&  Relative residuals \\
\hline
 4&1 &5 &0.172328\\ 
 \cline{2-4}
 &2&5&0.003307\\
  \cline{2-4}
 &3&1&3.165177E-06\\
  \cline{2-4}
 &4&1&2.408886e-10\\
\hline   
 16&1 &5 &  0.168451\\ 
 \cline{2-4}
 &2&5&0.002085\\
  \cline{2-4}
 &3&5&1.009380E-06\\
  \cline{2-4}
 &4&1&2.813398E-09\\
\hline 
\end{tabular}
\captionsetup{type=table}
\captionof{table}{Results for NNs training outside with $H/h=4$.}
\label{FAS_outside_H/h_4}
\end{minipage}

\subsubsection{Cost of training}
\begin{figure}[!ht]
\centering
\subfloat[Training and testing accuracy]{\includegraphics[scale=0.5]{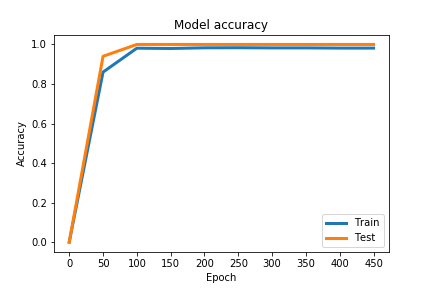}} 
\subfloat[Training and testing loss]{\includegraphics[scale=0.5]{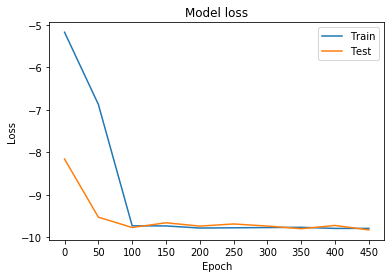}}
\caption{Plots of the accuracy at every 50 epochs and the log scale of loss values at every 50 epochs for one subdomain. We use $80\%$ of the number of vectors for training purpose and $20\%$ for testing. Here, we do the training inside the FAS. 
\label{figure: accuracy and loss for training inside}}
\end{figure}

\begin{figure}[!ht]
\centering
\subfloat[Training and testing accuracy]{\includegraphics[scale=0.5]{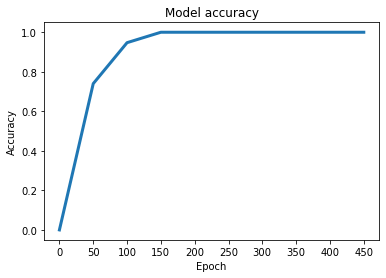}} 
\subfloat[Training and testing loss]{\includegraphics[scale=0.5]{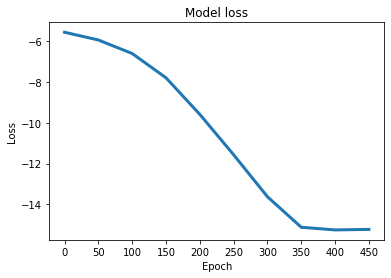}}
\caption{We do the training outside the FAS and use all the data/vectors for testing the accuracy as well as the loss values. We plot the accuracy at every 50 epochs and the log scale of loss values at every 50 epochs for one subdomain. \label{figure: accuracy and loss for training outside}}
\end{figure}

To get a sense of the cost, for the case of 4 subdomains and $H/h=2$, we display the accuracy and the loss values, characteristics provided by Keras, commonly used in training DNNs. Here, we only present the results for one of the four  subdomains since the outcomes are similar for all other  subdomains. Figure~\ref{figure: accuracy and loss for training inside} shows plots for training inside whereas Figure~\ref{figure: accuracy and loss for training outside} for the training outside FAS.

\newpage
\subsubsection{Results for different nonlinear coefficients $k=k(u)$}
Next, we consider 4 subdomains with $H/h=2$ for different coefficients $k$.  We use the same settings for the initial inputs and the neural networks as specified  at the beginning of this section. We see in Table~\ref{table1}, Table~\ref{table2}, Table~\ref{table3}, Table~\ref{table4}, Table~\ref{table5} and Table~\ref{table6} that the results are comparable 
for all of the used coefficients $k = k(u)$.
\vspace{0.5cm}

More specifically, the results for  $k=1+e^{-u}+x^2 +y^2$ are 
found in Table~\ref{table1}, Table~\ref{table2}, Table~\ref{table3}.\vspace{0.5cm}

\noindent\begin{minipage}{\linewidth}
\centering
\begin{tabular}{|c|c|c|}
\hline
FAS iteration& \# coarse iterations&  Relative residuals \\
\hline
1 & 5 & 0.128534    \\ \hline
2 & 5 & 0.000738     \\ \hline 
3 & 1 & 3.857236E-07   \\ \hline
\end{tabular}
\captionsetup{type=table}
\captionof{table}{Results for FAS using true operators}
\label{table1}
\end{minipage}

\noindent\begin{minipage}{\linewidth}
\centering
\begin{tabular}{|c|c|c|}
\hline
FAS iteration& \# coarse iterations&  Relative residuals \\
\hline
1 & 5 & 0.127007    \\ \hline
2 & 5 & 0.000706     \\ \hline
3 & 1 & 2.677881E-07    \\ \hline
\end{tabular}
\captionsetup{type=table}
\captionof{table}{Results for FAS with training inside}
\label{table2}
\end{minipage}

\noindent\begin{minipage}{\linewidth}
\centering
\begin{tabular}{|c|c|c|}
\hline
FAS iteration& \# coarse iterations& Relative residuals \\
\hline
1 & 5 & 0.141989   \\ \hline
2 & 5 & 0.001091   \\ \hline
3 & 5 & 7.363723E-07  \\ \hline
\end{tabular}
\captionsetup{type=table}
\captionof{table}{Results for FAS with training outside}
\label{table3}
\end{minipage}

Similarly, for $k=1+e^{-u}$, we have the results displayed in Table~\ref{table4}, Table~\ref{table5} and Table~\ref{table6}.
\vspace{1cm}

\noindent\begin{minipage}{\linewidth}
\centering
\begin{tabular}{|c|c|c|}
\hline
FAS iteration& \# coarse iterations&  Relative residuals \\
\hline
1 & 5 & 0.127675    \\ \hline
2 & 5 &  0.000701     \\ \hline
3 & 1 & 2.662172E-07     \\ \hline
\end{tabular}
\captionsetup{type=table}
\captionof{table}{Results for FAS with true operators}
\label{table4}
\end{minipage}

\noindent\begin{minipage}{\linewidth}
\centering
\begin{tabular}{|c|c|c|}
\hline
FAS iteration& \# coarse iterations&  Relative residuals \\
\hline
1 & 5 & 0.124039    \\ \hline
2 & 5 & 0.000535   \\ \hline
3 & 1 & 1.782604E-07     \\ \hline
\end{tabular}
\captionsetup{type=table}
\captionof{table}{Results for FAS with training inside}
\label{table5}
\end{minipage}

\noindent\begin{minipage}{\linewidth}
\centering
\begin{tabular}{|c|c|c|}
\hline
FAS iteration& \# coarse iterations&  Relative residuals \\
\hline
1 & 5 & 0.141077    \\ \hline
2 & 5 & 0.001084   \\ \hline
3 & 5 & 7.305067E-07     \\ \hline
\end{tabular}
\captionsetup{type=table}
\captionof{table}{Results for FAS with training outside}
\label{table6}
\end{minipage}

\subsubsection{Illustration of the computed approximate solutions}
We also provide illustration for the approximate solutions obtained from Newton method using the nonlinear operators from the training outside FAS and project them back to the fine level along with the true solutions in the following figures \ref{solution1}, \ref{solution2} and \ref{solution3} for several subdomain cases. Here, we use the true solution $u = x^2(1-x)^2 + y^2(1-y)^2$. These illustrations demonstrate the potential for using the trained DNNs as accurate discretization tool. 
\begin{center}
\begin{figure}[!ht]
\subfloat[True solution]{\includegraphics[scale=0.15]{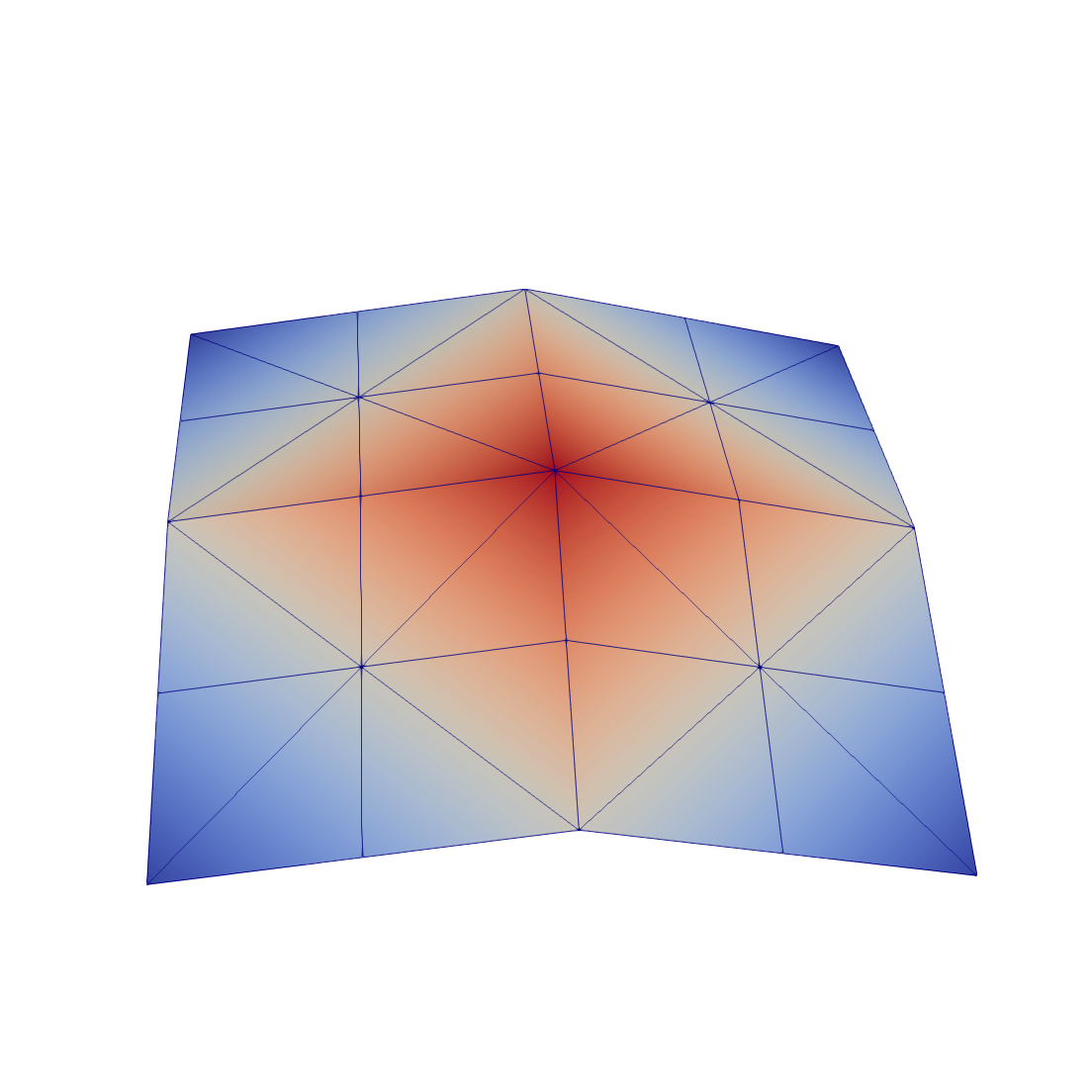}} 
\subfloat[Approximate solution]{\includegraphics[scale=0.15]{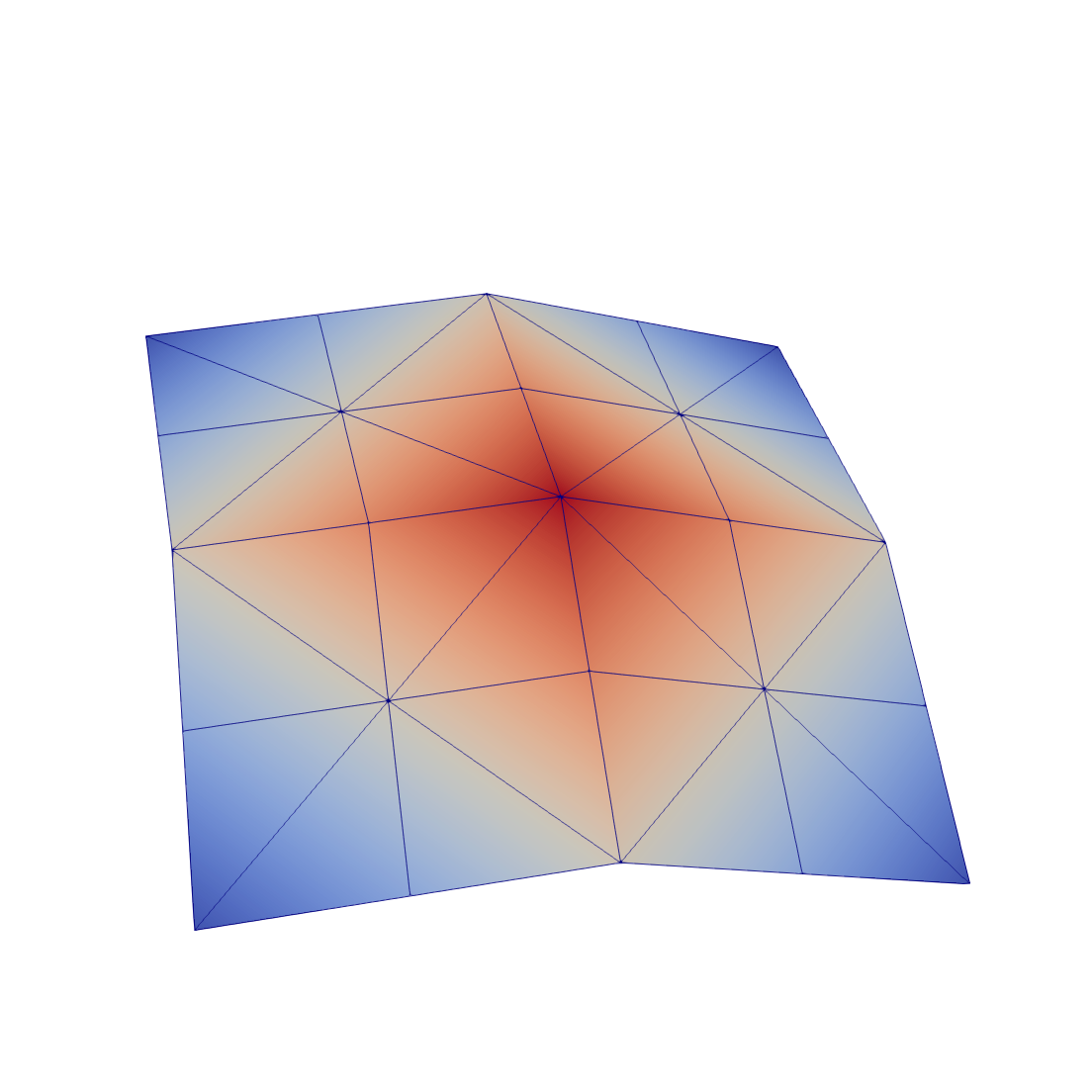}}
\caption{Plots of approximate and true solutions for 4 subdomains.}
\label{solution1}
\end{figure}
\end{center}

\begin{center}
\begin{figure}[!ht]
\subfloat[True solution]{\includegraphics[scale=0.15]{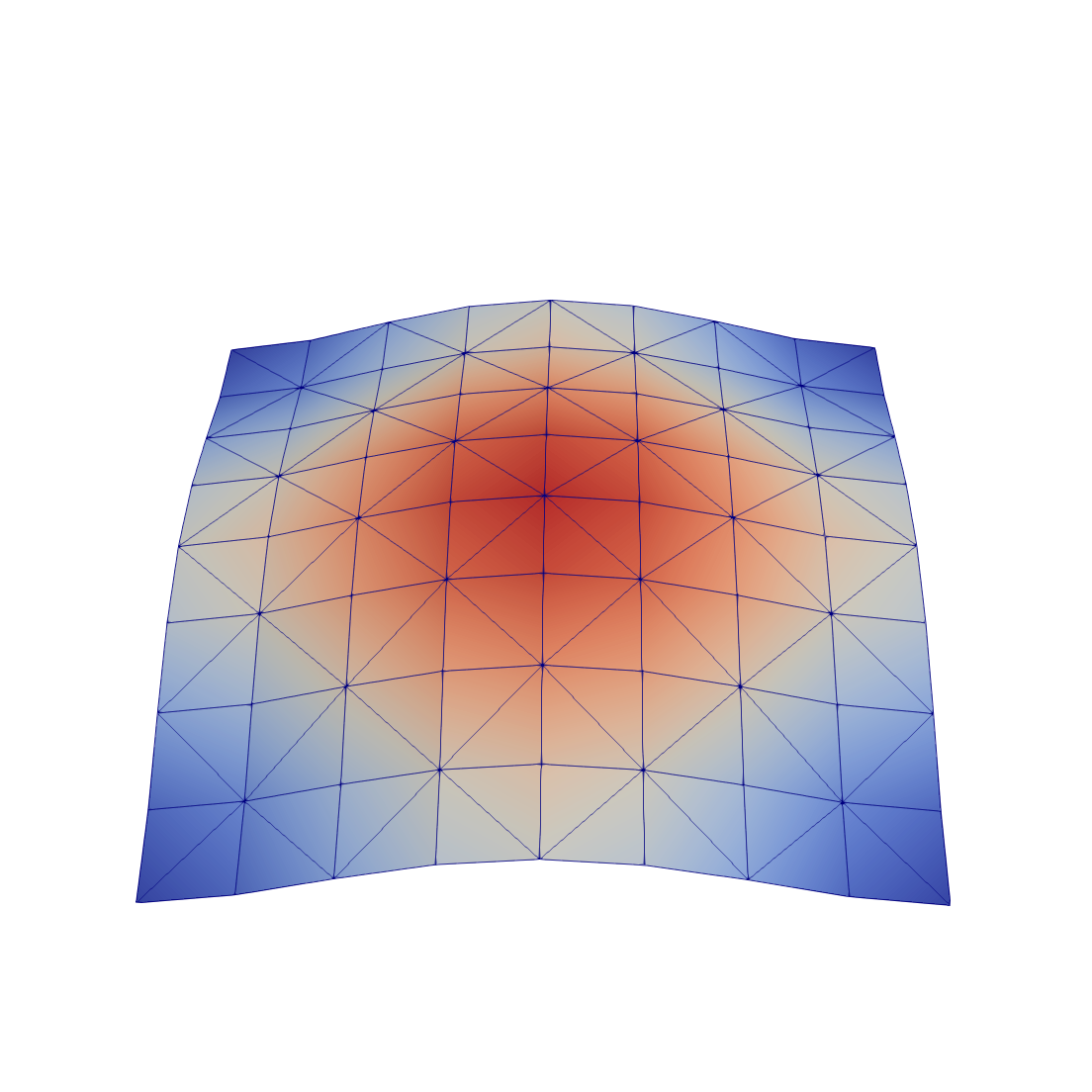}} 
\subfloat[Approximate solution]{\includegraphics[scale=0.15]{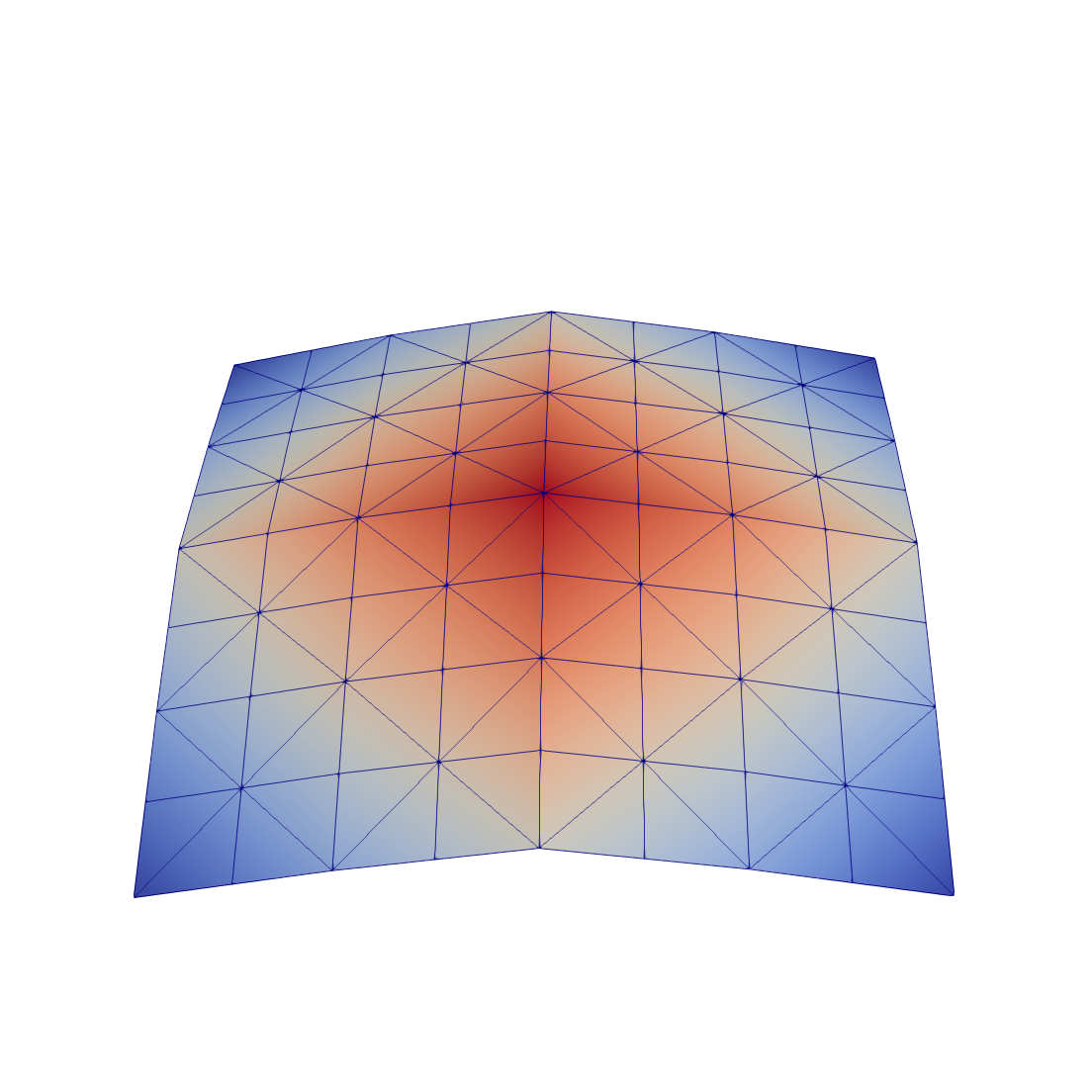}}
\caption{Plots of approximate and true solutions for 16 subdomains.}
\label{solution2}
\end{figure}
\end{center}

\begin{center}
\begin{figure}[!ht]
\subfloat[True solution]{\includegraphics[scale=0.15]{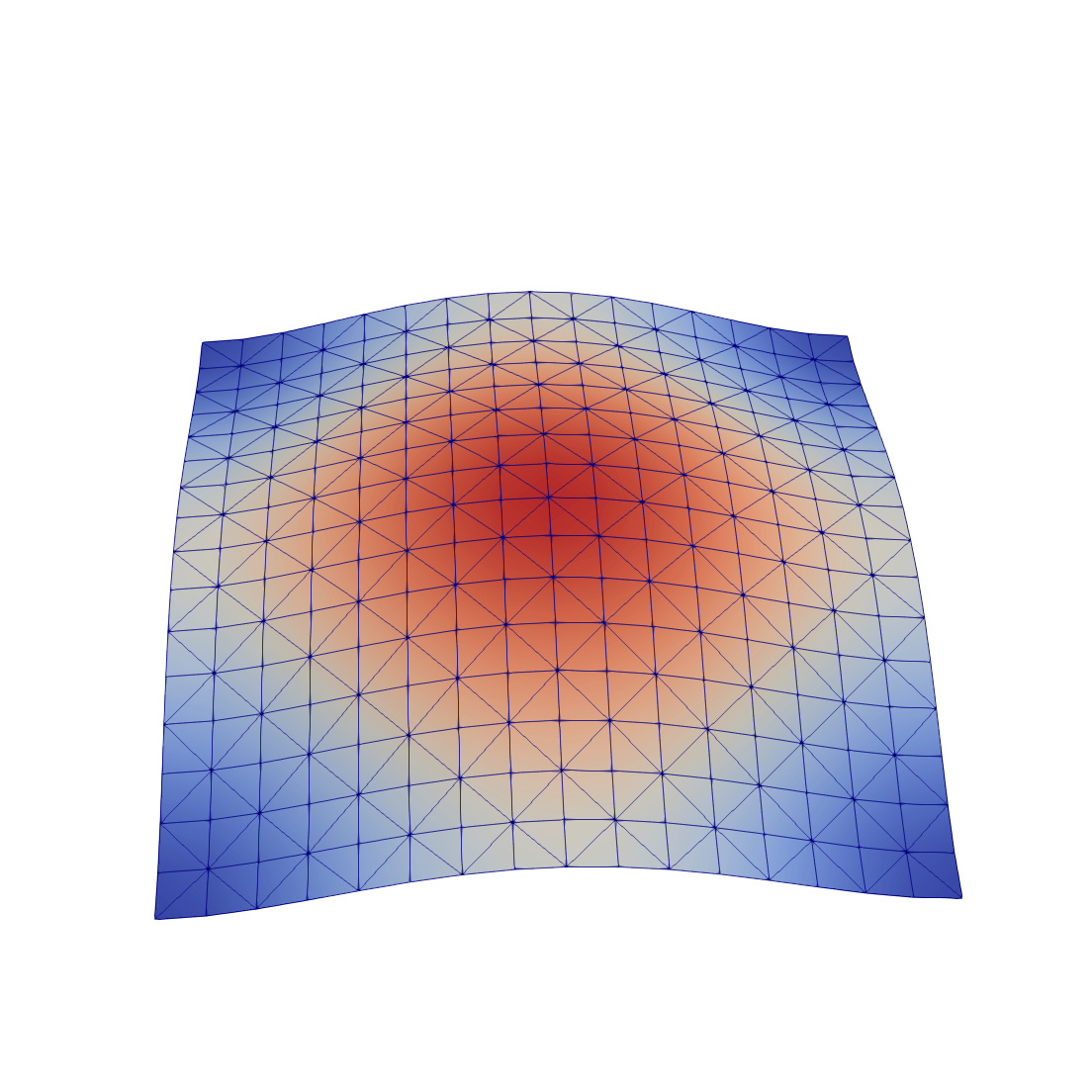}} 
\subfloat[Approximate solution]{\includegraphics[scale=0.15]{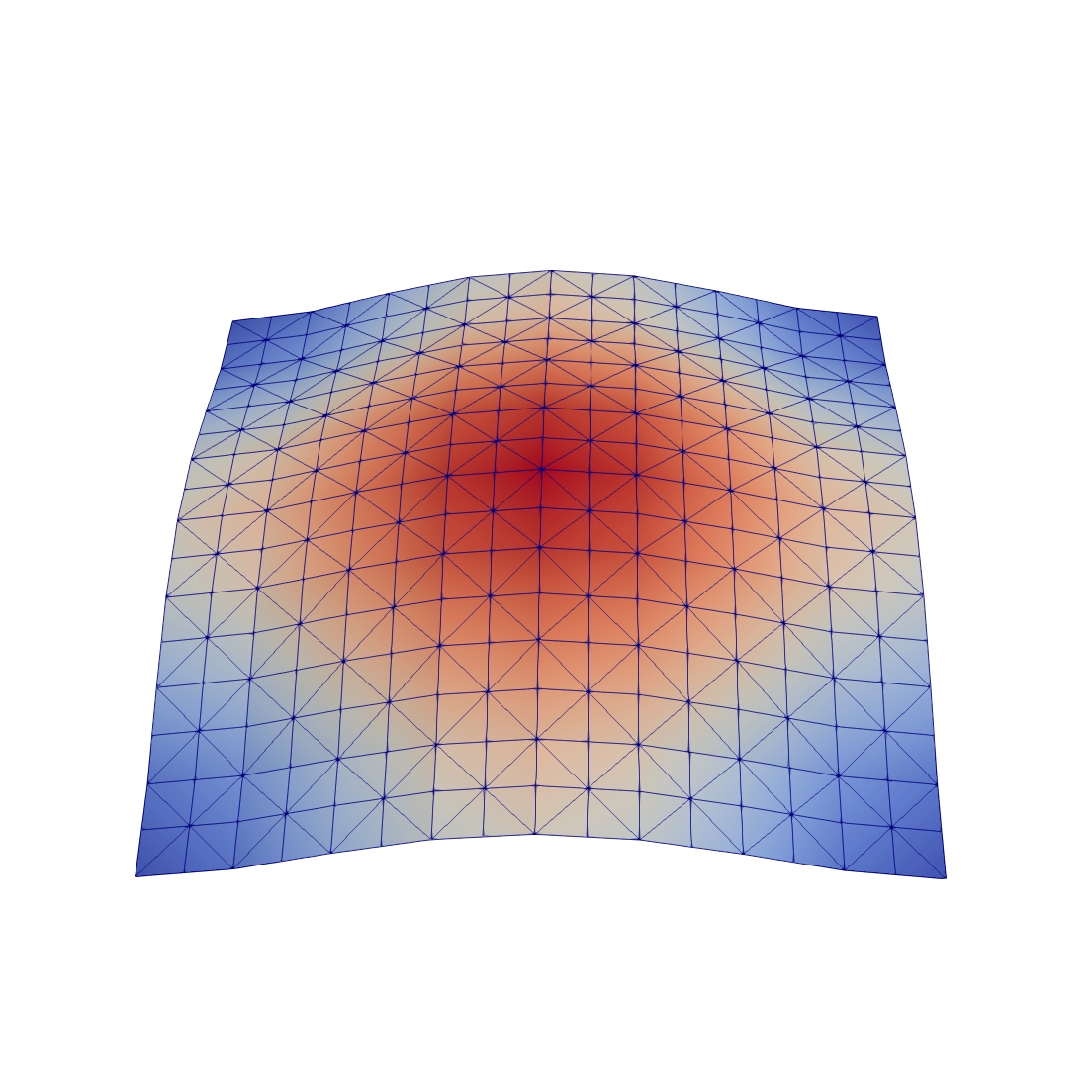}}
\caption{Plots of approximate and true solutions for 64 subdomains.}
\label{solution3}
\end{figure}
\end{center}

\begin{center}
\begin{figure}[!ht]
\subfloat[True solution]{\includegraphics[scale=0.3]{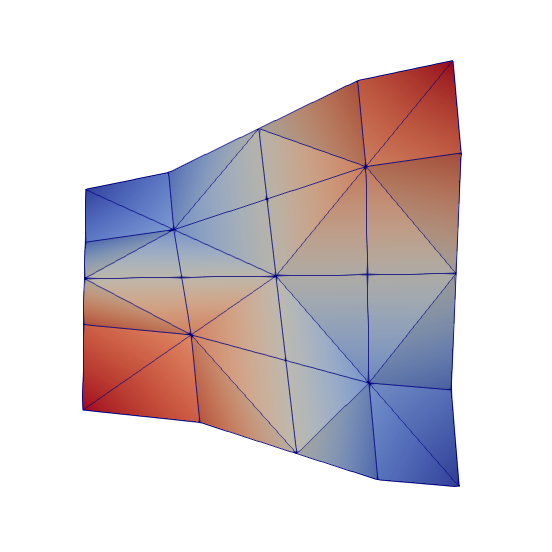}} 
\subfloat[Approximate solution]{\includegraphics[scale=0.33]{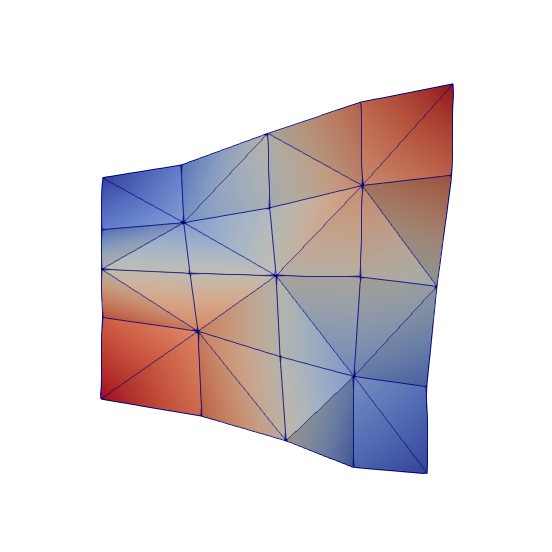}}
\caption{Plots of approximate and true solutions for 4 subdomains.}
\label{solution4}
\end{figure}
\end{center}

\begin{center}
\begin{figure}[!ht]
\subfloat[True solution]{\includegraphics[scale=0.3]{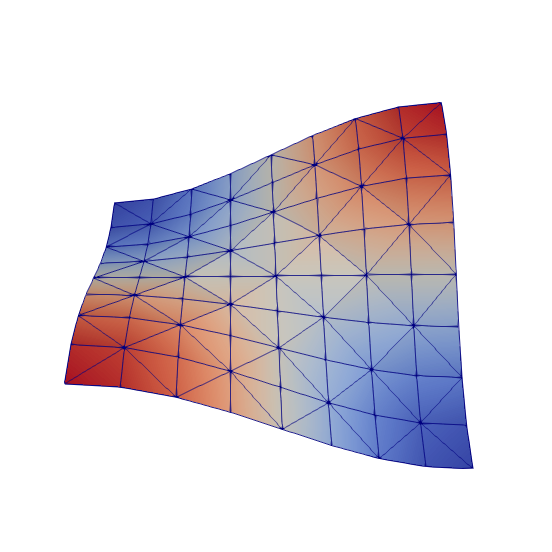}} 
\subfloat[Approximate solution]{\includegraphics[scale=0.35]{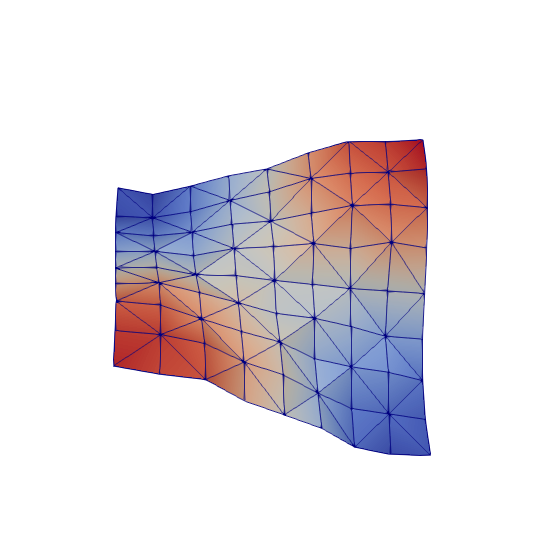}}
\caption{Plots of approximate and true solutions for 16 subdomains.}
\label{solution5}
\end{figure}
\end{center}

\begin{center}
\begin{figure}[!ht]
\subfloat[True solution]{\includegraphics[scale=0.28]{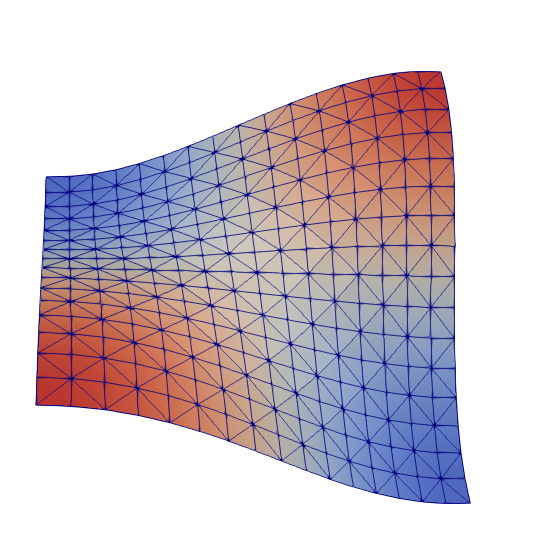}} 
\subfloat[Approximate solution]{\includegraphics[scale=0.28]{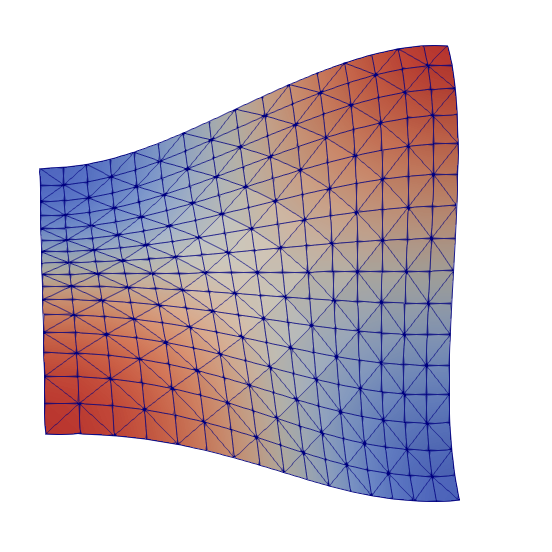}}
\caption{Plots of approximate and true solutions for 64 subdomains.}
\label{solution6}
\end{figure}
\end{center}

Similarly, we perform the test with the exact solution $u=\mbox{cos}(\pi x)\mbox{cos}(\pi y)$. Figure  \ref{solution4}, \ref{solution5}, and \ref{solution6} show the plots of approximate solutions and true solutions for different number of subdomains.

\newpage

\newpage
\section{Conlusions and future work} \label{conclusion}
The paper presented first encouraging results for approximating coarse finite element nonlinear operators for model diffusion reaction PDE in two dimensions.
These operators were successfully employed in a two-level FAS for solving the resulting system of nonlinear algebraic equations.
The resulting DNNs are quite expensive to replace the true Galerkin coarse nonlinear operators, however once constructed, one could in principle use them for solving the same type nonlinear PDEs with different r.h.s.
Upto a certain extent we can control the DNN complexity by choosing larger ratio $H/h$
and finally, it is clear that the training since it is local, subdomain-by-subdomain, and independent of each other, one can exploit parallelism in the training. Another viable option is to use convolutional DNNs instead of the currently employed fully connected ones. Also, a natural next step is to apply recursion, thus ending up with a hierarchy of trained coarse DNNs for use as coarse nonlinear discretization operators. 	
There is one more part where we can apply DNNs, namely to get approximations to the coarse Jacobians, $J_c(\bu_{c,T})$ (also done locally). Here, the input is the local coarse vector
$\bu_{c,\;T}$ and the output will be a local matrix $J_{c,T}(\bu_{c,T})$. 
It is also of interest to consider more general nonlinear, including stochastic, PDEs, which is a widely open area for future research. 

\section*{Disclaimer}
{\tiny 
This document was prepared as an account of work sponsored by an agency of the United States government. Neither the United States government nor Lawrence Livermore National Security, LLC, nor any of their employees makes any warranty, expressed or implied, or assumes any legal liability or responsibility for the accu- racy, completeness, or usefulness of any information, apparatus, product, or pro- cess disclosed, or represents that its use would not infringe privately owned rights. Reference herein to any specific commercial product, process, or service by trade name, trademark, manufacturer, or otherwise does not necessarily constitute or im- ply its endorsement, recommendation, or favoring by the United States government or Lawrence Livermore National Security, LLC. The views and opinions of authors expressed herein do not necessarily state or reflect those of the United States gov- ernment or Lawrence Livermore National Security, LLC, and shall not be used for advertising or product endorsement purposes.}

\end{document}